\documentclass[reqno,11pt]{amsart}

\usepackage[top=2.0cm,bottom=2.0cm,left=3cm,right=3cm]{geometry}
\usepackage{amsthm,amsmath,amssymb,dsfont}
\usepackage{mathrsfs,amsfonts,functan,extarrows,mathtools}
\usepackage[colorlinks]{hyperref}

\usepackage{marginnote}
\usepackage{xcolor}
\usepackage{stmaryrd}
\usepackage{esint}
\usepackage{graphicx}
\usepackage{bbm}
\usepackage{wasysym}
\usepackage{textcomp}
\usepackage{float}
\usepackage{wasysym}
\usepackage{cite}
\usepackage{tikz}
\usepackage{bookmark}
\usetikzlibrary{calc,decorations.pathreplacing}

\newtheorem{theorem}{Theorem}[section]
\newtheorem{proposition}{Proposition}[section]
\newtheorem{definition}{Definition}[section]
\newtheorem{corollary}{Corollary}[section]
\newtheorem{remark}{Remark}[section]
\newtheorem{lemma}{Lemma}[section]

\numberwithin{equation}{section}
\allowdisplaybreaks

\def\d{\mathrm{d}}

\def\R{\mathbb{R}}
\def\eps{\varepsilon}

\def\exp{\mathrm{exp}}

\def\L{\mathcal{L}}

\newcounter{wronumber}

\begin{document}
\title[Boundary layer equation]{Boltzmann boundary layer equation with Maxwell reflection boundary condition and applications to fluid limits}
\date{}
\author[Ling-Bing He]{Ling-Bing He}
\address[Ling-Bing He]{\newline Department of Mathematical Sciences, Tsinghua University, Beijing, 100084, P.R. China}
\email{hlb@tsinghua.edu.cn}

\author[Ning Jiang]{Ning Jiang}
\address[Ning Jiang]{\newline School of Mathematics and Statistics, Wuhan University, Wuhan, 430072, P. R. China}
\email{njiang@whu.edu.cn}

\author[Yulong Wu]{Yulong Wu${}^*$}
\address[Yulong Wu]{\newline School of Mathematics and Statistics, Wuhan University, Wuhan, 430072, P. R. China}
\email{yulong\string_wu@whu.edu.cn}

\thanks{${}^*$ Corresponding author \quad \today}

\maketitle

\begin{abstract}
  This paper investigates the Knudsen layer equation in half-space, arising from the hydrodynamic limit of the Boltzmann equation to fluid dynamics. We consider the Maxwell reflection boundary condition with accommodation coefficient $0<\alpha<1$. We restrict our attention to hard sphere collisions with angular cutoff, proving the existence, uniqueness, and asymptotic behavior of the solution in $L^{\infty}_{x,v}$. Additionally, we demonstrate the application of our theorem to the hydrodynamic limit through a specific example. In this expample, we derive the boundary conditions of the fluid equations using our theorem and the symmetric properties of the Knudsen layer equation for $\alpha\in(0,1]$ and $\alpha=O(1)$. These derivations differs significantly from the cases of specular and almost specular reflection. This explicitly characterizes the {\em vanishing sources set} defined in \cite{jiang2024knudsenboundarylayerequations}.\\


  \end{abstract}

\section{Introduction and main results}
\subsection{Problem setting}\label{Se1.1}
In the present paper, we study the following nonlinear Milne problem: the mass density $F$ is a function of $x\in \mathbb{R}_+$ and particle velocity $v=(v_1,v_2,v_3)\in \mathbb{R}^3$. Let $v_3$ stand for the velocity component along the $x$-axis. Then, $F(x,v)$ is governed by
  \begin{equation}\label{1.1.1}
      v_3 \partial_x F=Q(F,F)+\mathcal{S},
  \end{equation}
where $Q$ denotes the Boltzmann collision operator, defined by
    \begin{equation*}
      Q(F_1,F_2)=\int_{\mathbb{R}^3 \times \mathbb{S}^2}b(|v_\ast-v|,\omega )[F_1(v_\ast^{\prime})F_2(v^{\prime})-F_1(v_\ast)F_2(v)]dv_\ast d\omega.
      \end{equation*}
       Here $\omega \in \mathbb{S}^2$ is a unit vector and $v,v_\ast $ and $v',v'_\ast $ are the velocities of the particles before and after the collision, respectively, which satisfy the conservation of momentum and mass:
      $$\left\{
     \begin{aligned}
     &v'=v+[(v_\ast -v)\cdot w]w,\quad v_\ast'=v_\ast-[(v_\ast -v)\cdot w]w,\\
     &|v_\ast |^2+|v|^2=|v'_\ast |^2+|v'|^2.
     \end{aligned}\right.
     $$
     For simplicity, throughout this paper, we consider the hard sphere model and the Boltzmann collision kernel takes the form
     $$b(|v_\ast-v|,\omega )=|(v_\ast-v)\cdot\omega|.$$
Moreover, the equation \eqref{1.1.1} is subject to the Maxwell reflection boundary condition with the source term
 \begin{equation}\label{1.1.3}
      F(0,v) |_{v_3>0}=(1-\alpha)F(0,R_0 v)+\alpha \sqrt{2\pi} m(v) \int_{u_3<0} |u_3| F(0,u) \d u + \mathcal{R}(v),
  \end{equation}
and satisfies the far-field condition
\begin{equation}\label{1.1.6}
  F\to m \quad\text {as } x \to \infty,
\end{equation}
where $m$ denotes the  global Maxwellian
$$m(v)=\tfrac{1}{(2\pi)^{3/2}}e^{-\frac{|v|^2}{2}},$$
and $\mathcal{R}(v) $ is a given function satisfying if $\mathcal{R}(v)=0 $ for $v_3<0$. The reflection operator is defined as $R_0 v= v-2 (v\cdot n_0) n_0$, where $n_0=(0,0,-1)$ is the unit outer normal vector at $x=0$. The constant  $\alpha$ is the accommodation coefficient, which describes the portion of the molecules accommodate the state of the wall or specularly reflected. Here we assume that $0 < \alpha < 1$, since the endpoint cases $\alpha=0$ and $\alpha=1$ which correspond to specular reflection and completely diffusive reflection respectively have been treated in \cite{huang2023boundary} and \cite{huang2022boundary}.

Let $\Sigma:= \partial \R_+ \times \R^3 $ be the phase space boundary of $\R_+\times \R^3$. The phase boundary $\Sigma$ can be split by outgoing boundary $\Sigma_+$, incoming boundary $\Sigma_-$, and grazing boundary $\Sigma_0$:
\begin{equation*}
  \begin{aligned}
    \Sigma_+=\{(x,v):x \in \partial\R_+,v\cdot n_0>0\},\\
    \Sigma_-=\{(x,v):x \in \partial\R_+,v\cdot n_0<0\},\\
    \Sigma_0=\{(x,v):x \in \partial\R_+,v\cdot n_0=0\}.
  \end{aligned}
\end{equation*}
Let $\gamma_{\pm}F=\mathbf{1}_{\Sigma_{\pm}}F$. Linearizing around $m$ in the form $F=m+\sqrt{m}f$, then \eqref{1.1.1}, \eqref{1.1.3} and \eqref{1.1.6} can be rewritten as
\begin{equation}\label{1.1.7}
\begin{cases}
  v_3  \partial_x f +\mathcal{L}f=\Gamma (f,f)+S, \quad (x,v) \in \mathbb{R}_+\times \mathbb{R}^3,\\
  (L^R-\alpha L^D)f(0,v)=R(v), \\
  \lim_{x\to \infty} f(x,v) = 0,
  \end{cases}
  \end{equation}
with
\begin{equation*}
  S(x,v) =\tfrac{\mathcal{S}}{\sqrt{m}}, \quad R(v)=\tfrac{\mathcal{R}}{\sqrt{m}},
\end{equation*}
and
\begin{equation*}
  L^R=\gamma_--L\gamma_+, \quad L^D=L\gamma_+-K\gamma_+,
\end{equation*}
  where $L\gamma_+ f(0, v)=f(0 ,R_0 v)$, and 
  \begin{equation}\label{Kgamma}
    K\gamma_+ f (x,v)=\sqrt{2\pi m(v)} \int_{u_3<0}{|u_3| f(x,u )\sqrt{m(u)}}d u.
  \end{equation}
Moreover, the linearized Boltzmann operator $\mathcal{L}$ is defined by
\begin{equation*}
   \mathcal{L} f(v)=-\tfrac{1}{\sqrt{m}}\left(Q(m,\sqrt{m }f)+Q(\sqrt{m}f,m)\right)
   =\nu(v)f(v)-Kf(v).
  \end{equation*}
 Here $\nu(v)$ denotes the collision frequency, given by 
  \begin{equation}\label{coll-freq}
    \nu(|v|)=\int_{\mathbb{R}^3 \times \mathbb{S}^2 }|(v-u)\cdot \omega|m(v)d\omega dv_* \cong 1+|v|,
  \end{equation}
where $A \cong B$ implies $\tfrac{1}{C} B\leq A \leq C B$ for some harmless constant $C>0$. The operator $K f(v)$ can be decomposed into two parts:
\begin{equation*}
	\begin{aligned}
		K f (v) = - K_1 f (v) + K_2 f (v).
	\end{aligned}
\end{equation*}
Here the loss term $K_1 f (v)$ is
\begin{equation*}
	\begin{aligned}
		K_1 f(v) = m^\frac{1}{2} (v) \iint_{\mathbb{R}^3 \times \mathbb{S}^2} f (v_*) m^\frac{1}{2} (v_*)  |\omega \cdot (v_* - v)| d \omega d v_*  ,
	\end{aligned}
\end{equation*}
and the gain term $K_2 f (v)$ is
\begin{equation*}
	\begin{aligned}
		K_2 f(v) = & \iint_{\mathbb{R}^3 \times \mathbb{S}^2} \big[ m^{\frac{1}{2}} (v'_*) f (v') + m^{\frac{1}{2}} (v') f(v_*') \big] m^\frac{1}{2} (v_*) |\omega \cdot (v_* - v)| d \omega d v_*.
	\end{aligned}
\end{equation*}
The bilinear term is defined as
\begin{equation*}
  \Gamma(f,f)=\tfrac{1}{\sqrt{m}}Q(\sqrt{m}f,\sqrt{m}f).
\end{equation*}

It is well-known that the null space $\mathcal{N}$ of operator $ \mathcal{L}$ is generated by (see \cite{caflisch1980fluid}, for instance)
$$\chi_0 = \sqrt{m}, \quad \chi_j =v_j \sqrt{m}(j=1,2,3),\quad \chi_4= \tfrac{|v|^2-3}{2}\sqrt{m}.$$
We define the macroscopic projection $\mathbf{P} : L^2 \to \mathcal{N}$ by  
\begin{equation*}
  \mathbf{P}f=\sum_{j=0}^{4}\chi_j\int_{\mathbb{R}^3} f\chi_j d v.
\end{equation*}
For later use, we denote 
\begin{equation}\label{AB}
\begin{aligned}
A_{i j}(v) & =\left(v_{i} v_{j}-\tfrac{\delta_{i j}}{3}|v|^{2}\right) \sqrt{m(v)}, \\
B_{i}(v) & =\left(\tfrac{|v|^{2}-5}{2}\right) v_{i} \sqrt{m(v)}.
\end{aligned}
\end{equation}
It is obvious to know $A_{i j},  {B}_{i} \in \mathcal{N}^{\perp}$. There exist two scalar functions $a(|v|)$ and $b(|v|)$ such that (see \cite{desvillettes1994remark} for instance)
\begin{equation}\label{hatAB}
\begin{cases}
    \mathcal{L}^{-1}A=\hat{A}, \quad  \mathcal{L}^{-1}B=\hat{B};\\
    \hat{A}=\alpha(|v|)A, \quad \hat{B}=\beta(|v|)B,
\end{cases}
\end{equation}
where $\mathcal{L}^{-1}$ is the pseudo-inverse of $\mathcal{L}$ defined on  $\mathcal{N}^\perp$. We denote
  \begin{equation}\label{kappa}
    \begin{aligned}
    & \kappa_{1}:=\int_{\mathbb{R}^3}A _{i j}\mathcal{L}^{-1} A_{i j}dv>0,\quad i \neq j, \\
    & \kappa_{2}:=\int_{\mathbb{R}^3} B_{i} \mathcal{L}^{-1} B_{i}dv>0.
    \end{aligned}
  \end{equation}

\subsection{Main results}

\subsubsection{Linear problem}
Neglecting the nonlinear term in \eqref{1.1.7} leads the linearized problem:
\begin{equation}\label{ASS}
  \begin{cases}
  v_3\partial_xf+\mathcal{L}f=g,\quad (x,v)\in\mathbb{R}_+\times \mathbb{R}^3,\\
  (L^R-\alpha L^D)f(0,v)=r(v),  \\
  \lim_{x \to \infty}  f(x,v)=0.
  \end{cases}
  \end{equation}
Inspired by \cite{bardos1986milne,ukai2003nonlinear,coron1988classification,jiang2024knudsenboundarylayerequations}, this system is overdetermined and solvable only under specific solvability conditions. As stated in Lemma 2.1.1 of \cite{coron1988classification}, the following system is well-posed even without imposing the far-filed condition:
\begin{equation}\label{ASS1}
  \begin{cases}
  v_3\partial_x {f}+\mathcal{L} {f}=g,\quad (x,v)\in\mathbb{R}_+\times \mathbb{R}^3,\\
 (L^R-\alpha L^D)f(0,v)=r(v),   \\
 \int_{v_3<0}|v_3| f(0,v) \sqrt{m(v)} dv=\lambda,
  \end{cases}
  \end{equation}
for any constant $\lambda$. Moreover, there exists a unique $q^\infty(v) \in \mathcal{N}$, such that $\lim_{x\to \infty}f(x,v)=q^\infty(v)$. The value of $q^\infty(v)$ depends uniquely on $g(x,v)$, $r(v)$ and $\lambda$, indicating that the solution asymptotically approaches the null space $\mathcal{N}$ rather than zero at infinity. By subtracting $q^\infty(v)$ from the solution of the above system and  and using the fact that $(v_3\partial_x +\L) q^\infty(v) \equiv 0$, the system can be reformulated as:
  \begin{equation}\label{KL-eq}
    \begin{cases}
    v_3\partial_xf+\mathcal{L}f=g,\quad (x,v)\in\mathbb{R}_+\times \mathbb{R}^3,\\
    (L^R-\alpha L^D)f(0,v)=-(L^R-\alpha L^D)q^\infty(v)+r(v),  \\
    \lim_{x \to \infty}  f(x,v)=0.
    \end{cases}
    \end{equation}
   Here $q^\infty(v) \in \mathcal{N}$ is determined together with the solution of \eqref{KL-eq} and will be explicitly specified later. We will establish the well-posedness of \eqref{KL-eq} and prove that the solution to \eqref{KL-eq} does not depend on the choice of $\lambda$. For simplicity, the boundary term $\eqref{ASS1}_3$ is neglected here.
   
   We introduce the weight function
  \begin{equation}\label{wv}
    w(v)=\left<v\right>^\beta e^{\vartheta|v|^2},
  \end{equation}
  where $\left<v\right>=\sqrt{1+|v|^2}$ and $\beta,\vartheta$ are two given constants.
   We now state our main result:
  
\begin{theorem} \label{Thm-linear}
  Assume that $0<\alpha<1$, $\beta\ge 3$, and $0\leq \vartheta < \tfrac{1}{8}$. The source terms $g$ and $r$ are assumed to satisfy the solvability condition:
  \begin{equation}\label{1.1.8}
   g \in \mathcal{N}^{\bot}, \quad \text{and} \quad \!\int_{v_3>0}v_3\sqrt{m(v)}r(v)dv=0,
  \end{equation}
with
  \begin{equation}\label{1.1.25}
    \left\|e^{\sigma_0x}\nu^{-1}w g\right\|_{L^\infty_{x,v}}+|w r|_{L^\infty_{v}(\mathbb{R}^3_+)}<\infty.
  \end{equation}
 Then there exists a unique function $q^\infty(v)$ with the form
  \begin{equation}\label{1.1.19}
    \begin{aligned}
    q^\infty(v) &=\left({b_1^\infty v_1+b_2^\infty v_2+c^\infty(\tfrac{|v|^2-3}{2})}\right) \sqrt{m},\\
    |(b^\infty_1,b^\infty_2,c^\infty)| &\leq C  \left(\left\|e^{\sigma_0x}\nu^{-1} w  g\right\|_{L^\infty}+| w  r|_{L^\infty(\gamma_-)}\right),
  \end{aligned}
  \end{equation}
such that the Knudsen layer equation \eqref{KL-eq} admits a unique solution $f(x,v)$ satisfying
\begin{equation}\label{1.1.20}
  \left\|e^{\sigma x}w f\right\|_{L^\infty }+|w f(0)|_{L^\infty(\gamma)} \leq  C \left(\tfrac{1}{\sigma_0-\sigma}\left\|e^{\sigma_0x}\nu^{-1} w  g\right\|_{L^\infty}+|w r|_{L^\infty_{v}(\gamma_+)}\right),
\end{equation}
where $0<\sigma<\sigma_0$ is a constant, and $C >0$ is a constant depending only on $\alpha$ and $\sigma_0$.
\end{theorem}
\begin{remark}
   Theorem \ref{Thm-linear} explicitly characterizes the  vanishing sources set ($\mathbb{VSS}$) introduced in \cite{jiang2024knudsenboundarylayerequations}. Additionally, Theorem \ref{Thm-linear} can be applied to determine the boundary conditions of the fluid (and Prandtl-type) equations within the expansion methods (both Hilbert and Chapman-Enskog expansion) of the  Boltzmann equation. In Section \ref{Sec-application}, we provide an example of deriving the boundary conditions for the fluid equations in the hydrodynamic limit of the linearized Boltzmann equation.
\end{remark}
\begin{remark}
 We remark that $(L^R-\alpha L^D)a^\infty \equiv 0$ for any constant $a^\infty$, and thus the solution to the system \eqref{KL-eq} is independent of the choice of $a^\infty$. We emphasize that the uniqueness of $q^\infty(v)$ is understood in the sense of \eqref{1.1.19}. 
\end{remark}

\subsubsection{Nonlinear problem}
Based on the linear theory Theorem \ref{Thm-linear}, we now proceed to state the nonlinear theorem.
\begin{theorem}\label{Thm-nonlinear}
   Assume that $0<\alpha<1$, $\beta\geq 3$, and $0\leq \vartheta < \frac{1}{8}$. The source terms $S$ and $R$ are assumed to satisfy the solvability condition:
   \begin{equation}
    S \in \mathcal{N}^{\bot}, \quad \text{ and }\int_{v_3>0}v_3\sqrt{m(v)}R(v)dv=0.
   \end{equation}
Let $\sigma_0$ and $\delta_0$ be suitable small such that
  \begin{equation}
    \left\|e^{\sigma_0x}\nu^{-1}wS\right\|_{L^\infty}+|wR|_{L^\infty(\gamma_-)} \leq \delta_0.
  \end{equation}
  Then there exists a unique function $q^\infty(v)$ with the form
  \begin{equation}\label{1.1.13}
    \begin{aligned}
    q^\infty(v) &=\left(b_1^\infty v_1+b_2^\infty v_2+c^\infty (\tfrac{|v|^2-3}{2})\right) \sqrt{m},\\
    |(b^\infty_1,b^\infty_2,c^\infty)|  &\leq C  ( \left\|e^{\sigma_0x}\nu^{-1} w  S\right\|_{L^\infty }+| w  R|_{L^\infty (\gamma_-)} ),
  \end{aligned}
  \end{equation}
such that the nonlinear Knudsen layer problem
\begin{equation}\label{nonlinear-KL-eq}
  \begin{cases}
      v_3 \partial_x f +\mathcal{L}f=\Gamma (f,f)+S, \quad x \in \mathbb{R}_+\times \mathbb{R}^3,\\
    (L^R-\alpha L^D)f(0,v)|_{v_3>0}=-(L^R-\alpha L^D)q^\infty(v)+R(v),\\
    \lim_{x\to\infty}f = 0 ,
    \end{cases}
    \end{equation}
admits a unique solution $f(x,v)$. Moreover, $f(x,v)$ satisfies the following estimates
\begin{equation}\label{1.1.15}
  \left\|e^{\sigma x} w  f\right\|_{L^\infty}+|w f(0)|_{L^\infty(\gamma)} \leq  \tfrac{C}{ \sigma_0-\sigma }(\left\|e^{\sigma_0x}\nu^{-1} w  S\right\|_{L^\infty}+| w  R|_{L^\infty(\gamma_-)})<\infty,
\end{equation}
where $0<\sigma<\sigma_0$ is a constant, and $C >0$ is a constant depending only on $\alpha$ and $\sigma_0$.
\end{theorem}

\subsection{History and motivation} When considering the hydrodynamic limit of the Boltzmann equation or the zero viscousity limit of the Navier-Stokes equation in the domain with boundary, the boundary layer effect cannot be neglected, see \cite{jiang2021compressiblea,jiang2021compressibleb,Masmoudi-Jiang-2017-CPAM,Sammartino-Caflisch-1988-CMP1,Sammartino-Caflisch-1988-CMP2}. Significant progress has been made in the Knudsen layer problem for the Boltzmann equation over the past thirty years. For the incoming boundary condition, Bardos-Caflisch-Nicolaenko \cite{bardos1986milne} established the existence, uniqueness, and asymptotic behavior of the linear equation for hard sphere model, demonstrating exponential decay $e^{-cx}$ in the space $L^\infty(e^{cx};L^2(1+|v|)^{1/2})$. Later, Golse-Poupaud  \cite{golse1989stationary} extended these results to a more general collision kernel with $-2\leq\gamma<1$. Golse-Coron-Sulem \cite{coron1988classification} further investigated moving boundary for the hard sphere model, proving well-posedness and determining the number of the solvability conditions based on the Mach number of the far Maxwellian. Ukai-Yang-Yu \cite{ukai2003nonlinear} established the existence of nonlinear Knudsen boundary layers with the moving boundary for the hard sphere collision kernel for the nondegenerate cases, i.e. the Mach number $\neq 0,\pm1$. They introduced an artificial damping term $\xi_1 P^+_0\xi_1 f$, to achieve macroscopic dissipation, which can be removed under specific solvability conditions. Their results have also been extended to more general collision kernels; see \cite{Chen-Liu-Yang-2004-AA,Wang-Yang-Yang-2007-JMP,Yang-2011-JSP}. The nonlinear stability of the Knudsen layer equation has been studied in \cite{ukai2004nonlinear,deng2008pointwise,Wang-Yang-Yang-2006-JMP,Shota-2022-JDE}. Additionally, the theory of the boundary layer equations with phase transition has been developed in \cite{Liu-Yu-2013-ARMA} and \cite{Niclas-Golse-2021-ARMA}.

For Maxwell reflection boundary condition, Golse-Perthame-Sulem \cite{golse1988boundary} proved the existence, uniqueness, and asymptotic behavior in the space $L^2(e^{cx} dx ; L^2(1+|v|)^{1/2}) \cap L^\infty(e^{c x}dx;\\ L^2(dv))$ for specular reflection boundary condition, i.e. $\alpha=0$. Recently, Huang-Wang \cite{huang2022boundary} extended their result in the functional space $L^2_{x,v}\cap  L^\infty_{x,v}$. Golse-Coron-Sulem \cite{coron1988classification} proved the existence, uniqueness, and asymptotic behavior of linear Knudsen layer equation with a zero source term and Maxwell reflection boundary condition for $\alpha \in (0,1]$ in the functional space $L^\infty(e^{\gamma x}dx;L^2(|v_3|dv))$. Later, Huang-Wang \cite{huang2022boundary} extended the result to $L^2 _{x,v}\cap L^\infty_{x,v}$ with $\alpha=1$. In this paper, we fill the gap for $0<\alpha<1$ in $L^2 _{x,v}\cap L^\infty_{x,v}$ framework. Authors in \cite{jiang2024knudsenboundarylayerequations} studied the case $\alpha\in[0,1]$ and $-3<\gamma \leq 1$, proving the exponential decay $\exp\{-c x^{\frac{2}{3-\gamma}}-c|v|^2\}$ in $L^\infty_{x,v}$ framework. They considered the specific source terms, referred to as the vanishing sources set ($\mathbb{VSS}$). In our notations, the definition of $\mathbb{VSS}$ is defined as
\begin{equation*}
  \mathbb{VSS} =\{(g,r); g(x,v)\in \mathcal{N}^\perp \text{ and } r(v) \text{ in \eqref{ASS} such that } \lim_{x\to \infty}f(x,v)=0 \}.
\end{equation*}
In contrast, this paper, along with \cite{huang2022boundary,huang2023boundary}, focuses on general source terms. 

We now highlight the distinction between this paper and \cite{jiang2024knudsenboundarylayerequations}. In \cite{jiang2024knudsenboundarylayerequations}, the authors studied the system
  \begin{equation*}
  \begin{cases}
      v_3 \partial_x f +\mathcal{L}f=g, \quad (x,v) \in \mathbb{R}_+\times \mathbb{R}^3,\\
    (L^R-\alpha L^D)f(0,v)=r(v),\\
    \lim_{x\to \infty}f = 0.
    \end{cases}
    \end{equation*}
    As mentioned earlier, the far-field condition $\lim_{x\to +\infty}f =0$ is overdetermined for general source terms $(g,r)$. To address this issue, the authors introduced the so-called {\em vanishing sources set} ($\mathbb{VSS}$) and considered the sources $(g,r)\in \mathbb{VSS}$ to ensure that the solution vanishes at infinity. Similar approaches were employed in \cite{ukai2003nonlinear,Wang-Yang-Yang-2007-JMP,Chen-Liu-Yang-2004-AA}, where the incoming boundary condition was considered. Howover, in the fluid limits of the Boltzmann equation, how to determine the so-called $\mathbb{VSS}$ is unclear. This paper considers general source terms and provides an alternative description of the set $\mathbb{VSS}$, as stated in Theorem \ref{Thm-linear}.
    
    In Section \ref{Sec-application}, we apply Theorem \ref{Thm-linear} to the fluid limits of the Boltzmann equation. We slove the Knudsen layer equation and derive the boundary conditions of the fluid equations by using the symmetric properties of the linearized Boltzmann operator $\L$ and the boundary operator $L^R-\alpha L^D$. For simplicity, our example is on the linearized incompressible Navier-Stokes limit from the linearized Boltzmann equation. We note that this method is also applicable to other limits of the nonlinear Boltzmann equation with Maxwell reflection boundary condition, including the incompressible Euler and compressible Euler limits.

\subsection{Notations} We list the notations used in this paper. Throughout this paper, $C$ denotes a generic positive constant that may vary from line to line. We also introduce the following norms.
\begin{equation*}
  \begin{array}{l}
    \left\|f(x,v)\right\|_{L^2}=\left\|f(x,v)\right\|_{L^2_{x,v}}, \quad  \left\|f(x,v)\right\|_{L^\infty}=\left\|f\right\|_{L^\infty_{x,v}};\\
    \left\|f(x,v)\right\|_\nu=\left\|\sqrt{\nu}f\right\|_{L^2_{x,v}} \text{ or} \left\|\sqrt{\nu}f\right\|_{L^2_{v}}, 
  \end{array}
\end{equation*}
where $\|\cdot\|_{L^p_{x,v}}$ (respectively, $\|\cdot\|_{L^p_{v}}$) deontes the standard $L^p(\Omega \times \mathbb{R}^3)$ (respectively, $L^p( \mathbb{R}^3)$) norm.  For the phase boundary integration, $|f|_{L^2(\gamma_\pm )}$ denotes the $L^2(\gamma_\pm)$-norm and $|f|_{L^\infty(\gamma_\pm)}$ denotes the $L^\infty(\gamma_\pm)$-norm of $f$.

\subsection{Difficulties and Methodologies}
We outline the procedure for solving the Knudsen layer equation \eqref{KL-eq}. Noting that $q^\infty(v) \in \mathcal{L}$, it follows that $(v_3\partial_x + \L) q^\infty \equiv 0$. Consequently, \eqref{KL-eq} is equivalent to the following system
\begin{equation}\label{KL-infty}
  \begin{cases}
  v_3\partial_x {f}+\mathcal{L} {f}=g,\quad (x,v)\in\mathbb{R}_+\times \mathbb{R}^3,\\
 (L^R-\alpha L^D)f(0,v)=r(v),   \\
\lim_{x\to \infty} f(x,v)= q^\infty(v).
  \end{cases}
  \end{equation}
While investigating the problem \eqref{KL-infty} over half-space $x\in(0,\infty)$, we will consider \eqref{KL-infty} in a finite slab $x\in(0,d)$ to derive the $L^\infty_{x,v} $ estimates using the method of characteristics. However, since the far-filed condition is unknown, the natural boundary condition at $x=d$ is the specular reflection boundary. At $x=0$, the boundary condition is mixed: some particles undergo specular reflection, while others experience diffuse reflection. This creates more complex trajectories compared to the endpoint cases ($\alpha=0,1$), as discussed in \cite{huang2022boundary,huang2023boundary}, making it challenging to derive the $L^\infty_{x,v}$ estimates by the trajectory method. As inspired by \cite{coron1988classification}, we initially work on the following auxiliary system with the prescribed mass flux condition. In this framework, the diffusive reflection part is treated as a source term, which helps address the difficulties in obtaining  $L^\infty_{x,v}$ estimates. Finally, we will prove Theorem \ref{Thm-linear} and obtain the uniform estimates using a uniqueness arguement. Specifically, we define
\begin{equation}\label{1.1.22}
  P_\gamma f(0,v):=\int_{v_3<0}|v_3|f(0,v)\sqrt{m}dv=\lambda,
\end{equation}
for given $\lambda \in \R$ and consider the following auxiliary system, with prescribed mass flux at $x=0$:
\begin{equation}\label{Auxi-eq}
  \begin{cases}
  v_3\partial_x {f}+\mathcal{L} {f}=g,\quad (x,v)\in\mathbb{R}_+\times \mathbb{R}^3,\\
 (L^R-\alpha L^D)f(0,v)=r(v),  \\
P_\gamma f(0,v)=\lambda.
  \end{cases}
  \end{equation}
We emphasize that these difficulties do not arise under pure specular or pure diffuse boundary conditions, which have been studied in \cite{huang2023boundary} and \cite{huang2022boundary}. In these cases, the trajectories are much simpler and can be handled without considering the auxiliary system. Additionally, in \cite{jiang2024knudsenboundarylayerequations}, because the authors restrict the sources in the set $\mathbb{VSS}$ (admissible source set), the solution naturally vanishes at infinity and therefore can be approximated by incoming boundary condition at $x=d$. To characterize the set $\mathbb{VSS}$, we consider general sources here. 

Third, to construct the solution to the auxiliary system in the $L_{x,v}^\infty$ framework, we analyze the problem in a finite slab with specular reflection boundary condition at $x =d$. Specifically, we consider
\begin{equation}\label{Ap-eq-d}
  \begin{cases}
   v_3\partial_xf+ \mathcal{L}f=g,\quad (x,v)\in (0,d) \times \mathbb{R}^3,\\
   (L^R-\alpha L^D)f(0,v)=r(v), \\
   f(d,v)|_{v_3<0}=L\gamma_+ f(d,v),\\
 P_\gamma f(0,v)=\lambda,
  \end{cases}
  \end{equation}
where we use the notation $L\gamma_+ f(d,v)= f(d, R_d v)$, with $R_d v=v-2(v\cdot n_d)n_d$, and $n_d=(0,0,1)$ denoting the unit outer normal vector at $x=d$.
 
Fourth, we consider the approximate problem with an artificial damping term $\epsilon f$ to establish macroscopic coercivity. The term will be removed by taking  $\epsilon\to 0$ using compactness arguments.
\begin{equation}\label{Ap-eq-eps}
  \begin{cases}
   \epsilon f^\epsilon + v_3\partial_x f^\epsilon  + \mathcal{L}f^\epsilon =g,\quad (x,v)\in (0,d) \times \mathbb{R}^3,\\
   f^\epsilon(0,v) |_{v_3>0}=(1-\alpha)L\gamma_+ f^\epsilon(0,v)+ \alpha\lambda\sqrt{2\pi m(v)}+r(v), \\
f^\epsilon(d,v)|_{v_3<0}=L\gamma_+ f^\epsilon(d,v),\\
  \end{cases}
  \end{equation}

Fifth, we consider the approximate problem with an artificial damping and the parameter $a \in [0,1]$. We first construct the solution of \eqref{Con-eq} with $a=0$ and then extend the solution to $a=1$ by the contraction arguments.
\begin{equation}\label{Con-eq}
  \begin{cases}
   \epsilon f^\epsilon + v_3\partial_x f^\epsilon  + \nu f^\epsilon =a Kf^\epsilon + g,\quad (x,v)\in (0,d) \times \mathbb{R}^3,\\
   f^\epsilon(0,v) |_{v_3>0}=(1-\alpha)L\gamma_+ f^\epsilon(0,v)+ \alpha\lambda\sqrt{2\pi m(v)}+r(v), \\
 f^\epsilon(d,v)|_{v_3<0}=L\gamma_+ f^\epsilon(d,v).
  \end{cases}
  \end{equation}

Based on Step 2 to 5, we can construct the solution to \eqref{Auxi-eq}. More precisely, we have the following lemma.
\begin{lemma}\label{Lmm-auxi} Under the same assumption in Theorem \ref{Thm-linear}, for any $\lambda \in \R$, the auxiliary system \eqref{Auxi-eq} admits a unique solution $f(x,v)$. Moreover, there exists a unique $q^\infty\in \mathcal{N}$ satisfies
\begin{equation}
q^\infty=\left(a^\infty+b^\infty_1v_1+b^\infty_2v_2+c^\infty(\tfrac{|v|^2-3}{2})\right)\sqrt{m},
\end{equation}
where $(a^\infty,b^\infty_1,b^\infty_2,c^\infty)$ are  constants uniquely determined by $g,r(v)$ and $\lambda$, satisfying
\begin{equation}\label{1.26}
  \left\|e^{\sigma x} w  (f-q^\infty)\right\|_{L^\infty }+| w  (f(0)-q^\infty)|_{L^\infty(\gamma)} \leq C \left(\tfrac{1}{\sigma_0-\sigma}\left\|e^{\sigma_0x}\nu^{-1} w  g\right\|_{L^\infty }+|w r|_{L^\infty_{v}(\gamma +)}+|\lambda|\right),
\end{equation}
\begin{equation}\label{1.27}
  |(a^\infty,b^\infty_1,b^\infty_2,c^\infty)| \leq C\left(\left\|e^{\sigma_0x}\nu^{-1} w  g\right\|_{L^\infty}+| w  r|_{L^\infty(\gamma_-)}+|\lambda|\right),
\end{equation}
where $0<\sigma<\sigma_0$ is a constant, and $C >0$ is a constant depending only on $\alpha$ and $\sigma_0$.
\end{lemma}
Note that the estimates in Lemma \ref{Lmm-auxi} depend on the choice of $\lambda$. The final step is to construct the solution of \eqref{KL-eq} and obtain the uniform in $\lambda$ estimates. This is achieved through a uniqueness arguement, as discussed in subsection \ref{Sec-linear-KL}. Thus, we conclude Theorem \ref{Thm-linear}.

 In Section \ref{Sec-application}, we apply Theorem \ref{Thm-linear} to the fluid limits of the Boltzmann equation, illustrating how to derive the boundary conditions of the fluid equations by solving the Knudsen layer problem. The derivation of the boundary conditions essencially relies on the symmetric properties of the linearized Boltzmann operator $\L$ and the boundary operator $L^R-\alpha L^D$.

\section{Existence of the auxiliary system \eqref{Auxi-eq}}
In this section, we establish the existence, uniqueness, and asymptotic behavior of the auxiliary system \eqref{Auxi-eq}, i.e.,
\begin{equation*}
  \begin{cases}
   v_3\partial_xf+ \mathcal{L}f=g,\quad (x,v)\in \mathbb{R}_+ \times \mathbb{R}^3,\\
(L^R-\alpha L^D)f(0,v)=r(v),\\
  \int_{v_3<0}|v_3|f(0,v_3)\sqrt{m}dv=\lambda.
  \end{cases}
  \end{equation*}
The proof is divided into the following subsections.

\subsection{Existence of the Approxiamte system \eqref{Ap-eq-eps}}
In this subsection, we establish the existence of the approximate system \eqref{Ap-eq-eps}. We begin by solving the boundary value problem \eqref{Con-eq}, which depends on the parameter $a \in [0,1]$ and is equivalent to
\begin{equation}\label{Con-eq1}
  \begin{cases}
\mathscr{L}_a f^\epsilon:=\epsilon f^\epsilon + v_3\partial_x f^\epsilon+\nu(v)f^\epsilon= a K f^\epsilon+g,\\
  f^\epsilon(0,v)|_{v_3>0}=(1-\alpha)L\gamma_+f^\epsilon(0,v) +\tilde{r}(v),\\
f^\epsilon(d,v)|_{v_3<0}=L\gamma_+f^\epsilon(d,v),
  \end{cases}
  \end{equation}
where $\epsilon \in (0,1]$ and 
\begin{equation}\label{tilde-r}
  \tilde{r}(v)=r(v)+\alpha\lambda\sqrt{2\pi m(v)}.
\end{equation}
We define $\Omega_d:=(0,d)$ with $d \geq 1$. To derive the $L^\infty_{x,v}$ estimates, we define $ h^\epsilon(x,v)=wf^\epsilon(x,v)$, where $w$ is defined in \eqref{wv}. Consequently, equation \eqref{Con-eq1} can then be rewritten as
\begin{equation}\label{Con-eq2}
  \begin{cases}
  \epsilon h^\epsilon + v_3\partial_xh^\epsilon+\nu(v)h^\epsilon= a K_wh^\epsilon+wg,\quad (x,v)\in(0,d) \times \mathbb{R}^3,\\
  h^\epsilon(0,v)|_{v_3>0}=(1-\alpha)L\gamma_+h^\epsilon(0,v)+w\tilde{r}(v),\\
h^\epsilon(d,v)|_{v_3<0}=L\gamma_+ h^\epsilon(d,v),
  \end{cases}
  \end{equation}
where $K_w h=w K(\tfrac{h}{w})$.

To express the solution along characteristics, we introduce the standard {\em backward characteristics} defined as follows:
\begin{definition}
  For any given $[t,x,v]$, let $[X(s),V(s)]$ be the {\em backward characteristics} for \eqref{Con-eq2}, which are determined by 
  \begin{equation*}
    \begin{cases}
      \frac{dX(s)}{ds}=V_3(s), \quad \frac{dV(s)}{ds}=0,\\
      [X(t),V(t)]=(x,v).
    \end{cases}
  \end{equation*}
\end{definition}
The solution is given by
\begin{equation*}
  [X(s),V(s)]=[X(s;t,x,v),V(s;t,x,v)]=[x-(t-s)v_3,v].
\end{equation*}
Now, for each $(x,v)$ with $x\in \overline{\Omega}_d$ and $v_3\neq 0$, we define the {\em backward exit time} $t_b$ as the last moment when the backward trajectory $[X(-\tau;0,x,v ),V(-\tau;0,x,v)]$ remains within $\overline{\Omega}_d$:
\begin{equation*}
  t_b(x,v)=\sup\{s>0:x-\tau v_3 \in \overline{\Omega}_d \text{ for } 0 \leq \tau \leq s\}.
\end{equation*} 
 We also define 
 \begin{equation*}
  x_b(x,v)=x(t_b)=x-t_b(x,v)v_3 \in \Omega_d.
 \end{equation*}
For $x\in \overline{\Omega}_d, (x,v) \notin \gamma_0 \bigcup \gamma_-$ and $(t_0,x_0,v_0)=(t,x,v)$, we inductively define 
\begin{equation*}
  (t_{k+1},x_{k+1},v_{k+1})=(t_k-t_b(x_k,v_k),x_b(x_k,v_k),R_{x_{k+1}}v_k) \quad  \text{for} \quad \! k \geq 1.
\end{equation*} 
The {\em back-time cycle} is defined as follows:
\begin{equation*}
  \begin{cases}
    X_{cl}(s;t,x,v)=\sum_{k\geq 1}\mathbf{1}_{[t_{k+1}t_k)}(s)\{x_k-v_{k,3}(t_k-s)\};\\
    V_{cl}(s;t,x,v)=\sum_{k \geq 1}\mathbf{1}_{[t_{k+1}t_k)}(s)v_k.
  \end{cases}
\end{equation*}
It is easy to see that for $k\geq 1$ and $(x,v) \notin \gamma_0 \bigcup \gamma_-,$ it holds that
\begin{equation}\label{2.0.7}
  \begin{array}{l}
    x_k=\frac{1-(-1)^k}{2}x_1+\frac{1+(-1)^k}{2}x_2;\\
    t_k-t_{k+1}=t_1-t_2=\frac{d}{|v_{0,3}|}> 0.
  \end{array}
\end{equation}

We now present the following standard $L^\infty_{x,v}$ a priori estimates, obtained using the method of characteristics.
\begin{lemma}[A priori $L^\infty_{x,v}$ estimates]\label{Lmm-apriori-infty}
  Let $\epsilon \in [0,1]$, $a \in[0,1]$, $\alpha\in(0,1)$, and $h^\epsilon$ be the solution to \eqref{Con-eq2}. Then, the following inequality holds
  \begin{equation}\label{apriori-infty}
      \left\|h^\epsilon\right\|_{L^\infty} +|h^\epsilon|_{L^\infty (\gamma_+)}\leq C\left(\left\|\nu^{-1}wg\right\|_{L^\infty }+|w \tilde{r}|_{L^\infty (\gamma_-)}+\left\|\tfrac{h}{w}\right\|_{L^2} \right),
  \end{equation}
  where $C>0$ is a constant independent of $\epsilon$ and $a$.
\end{lemma}
\begin{proof}
The proof colsely follows Lemma 3.3 in \cite{huang2023boundary}, since the boundary conditions of \eqref{Con-eq2} contain only the specular reflection part. We omit the details here for brevity.
\end{proof}
\begin{lemma}[A priori $L^2_{x,v}$ estimates]\label{Lmm-apriori-L2}
  Let $\epsilon \in (0,1]$, $a \in[0,1]$, $\alpha\in(0,1)$, and  $f^\epsilon$ be the solution of \eqref{Con-eq1}
  with $ \left\|wf^\epsilon\right\|_{L^\infty}+|wf^\epsilon|_{L^\infty(\gamma)}<\infty.$ We have
  \begin{equation}\label{apriori-L2}
    \alpha |f^\epsilon|^2_{L^2(\gamma_+)} +\left\|f^\epsilon\right\|^2_{L^2}\leq \tfrac{C_{\epsilon}}{\alpha}\left( \|g \|^2_{L^2}+|\tilde{r}|^2_{L^2(\gamma_-)}\right),
  \end{equation}
  where the constant $C_{\epsilon}>0$ depends on $\epsilon$.
\end{lemma}
\begin{proof}
  Multiplying \eqref{Con-eq1} by $f^\epsilon$ and integrating over $x$ and $v$, one obtains that
  \begin{multline*}
      \epsilon\left\|f^\epsilon\right\|^2_{L^2}+\iint_{\Omega_d \times \R^3}\partial_x(\tfrac{1}{2} v_3 {f^\epsilon}^2) dx dv  + a \iint_{\Omega_d \times \R^3} (\nu  {f^\epsilon}- K f^\epsilon)f^\epsilon  dx dv \\
       + (1-a) \left\|f^\epsilon\right\|^2_\nu  = \iint_{\Omega_d \times \R^3}f^\epsilon g dx dv.
      \end{multline*}
  By the boundary condition in \eqref{Con-eq1}, it infers that
  \begin{equation*}
    \begin{aligned}
      \iint_{\Omega_d \times \R^3}\partial_x(\tfrac{1}{2} v_3 {f^\epsilon}^2) dx dv
      = &\tfrac{1}{2}\int_{ \R^3} v_3 {f^\epsilon}^2(d,v)  dv - \tfrac{1}{2}\int_{ \R^3} v_3 {f^\epsilon}^2(0,v)  dv\\
      = &-\int_{v_3<0} v_3 {f^\epsilon}^2 (0,v) dv -\int_{v_3>0} v_3 \big((1-\alpha){f^\epsilon}(R_0,v)+ \tilde{r}(v)\big)^2 dv\\
      \geq & \alpha |f^\epsilon|^2_{L^2(\gamma_+)} +\tfrac{C}{\alpha}|\tilde{r}|^2_{L^2(\gamma_-)},
    \end{aligned}
  \end{equation*}
  where we have used the fact that $\int_{ \R^3} v_3 {f^\epsilon}^2(d,v) dv=0$ by oddness.

  The coercivity of the linear Boltzmann operator $\L$ implies that
  \begin{equation}\label{2.2.13}
    \iint_{\Omega_d \times \R^3} f^\epsilon \L f^\epsilon  dxdv  \geq c_0 \left\|(\mathbf{I-P})f^\epsilon\right\|^2_\nu,
  \end{equation}
  for $c_0>0$. Moreover, H{\"o}lder's inequality shows 
  \begin{equation*}
    \iint_{\Omega_d \times \R^3}f^\epsilon g dx dv \leq \tfrac{\epsilon}{4}\|f^\epsilon\|_{L^2}^2 + C_\epsilon \|g\|_{L^2}^2.
  \end{equation*}
  As a result, one obtains that
  \begin{equation}\label{apriori-L21}
    \left\|\mathscr{L}_a^{-1}g\right\|^2_{L^2}=\left\|f^\epsilon\right\|^2_{L^2}\leq \tfrac{C_\epsilon}{\alpha}\left(\left\|g\right\|^2_{L^2}+|\tilde{r}|^2_{L^2(\gamma_-)}\right).
  \end{equation}
  The proof of Lemma \ref{Lmm-apriori-L2} is complete.
\end{proof}
Next, we will use an iterative method to construct the solution of \eqref{Con-eq1} for $a\in[0,1]$.
\begin{lemma}\label{Lmm-Con-eq}
  Let $\epsilon\in (0,1]$, $d\geq 1$ and $\alpha\in(0,1)$. Assume ${\left\|\nu^{-1}wg\right\|_{L^\infty }+|w \tilde{r}|_{L^\infty (\gamma_-)}} < \infty.$ Then there exists a unique solution to \eqref{Con-eq1} for all $a\in[0,1]$, satisfying
  \begin{equation}\label{2.2.1}
    \left\|wf^\epsilon\right\|_{L^\infty }+|wf^\epsilon|_{L^\infty(\gamma)} \leq \tfrac{C_{\epsilon,d}}{\alpha}\{\left\|\nu^{-1}wg\right\|_{L^\infty }+|w \tilde{r}|_{L^\infty (\gamma_-)}\},
  \end{equation}
where the constant $C_{\epsilon,d}>0$ depends only on $\epsilon $ and $d$.
\end{lemma}
\begin{proof}
  We divide our proof into the following steps:
  
  \textbf{Step 1. Existence of $\mathscr{L}^{-1}_0$.}
  
  First, we consider the following iterative system:
  \begin{equation}\label{2.2.3}
    \begin{cases}
      \epsilon f^{n+1} + v_3\partial_xf^{n+1}+\nu(v)f^{n+1}=g,\\
    f^{n+1}(0,v)|_{v_3>0}=(1-\alpha)L\gamma_+f^n(0,v)+\tilde{r}(v),\\
  f^{n+1}(d,v)|_{v_3<0}=L\gamma_+f^{n+1}(d,v),
    \end{cases}
    \end{equation}
  where $n \in \mathbb{N}$ and we set $f_0 \equiv 0$. We solve \eqref{2.2.3} inductively using the method of characteristics. Define $h^{n}=w(v)f^n$. The system for $h^{n+1}$ is given by:
   \begin{equation}\label{2.2.4}
    \begin{cases}
    \epsilon h^{n+1} + v_3\partial_xh^{n+1}+\nu(v)h^{n+1}=wg,\\
    h^{n+1}(0,v)|_{v_3>0}=(1-\alpha)L\gamma_+h^{n}(0,v)+w\tilde{r}(v),\\
  h^{n+1}(d,v)|_{v_3<0}=L\gamma_+h^{n+1}(d,v).
    \end{cases}
    \end{equation}
    We write the solution of \eqref{2.2.4} along the characteristics by Duhamel's principle as
   \begin{equation}\label{2.2.5}
    \begin{aligned}
      h^{n+1}(x,v)&=(1-\alpha)e^{-\nu_\epsilon(v)(t-t_2)}h^{n}(0,v)+e^{-\nu_\epsilon(v)(t-t_2)}\tilde{r}(R_0 v)\\
      &+\int_{t_2}^{t_1}e^{-(t-s)\nu_\epsilon}(wg)(X_{cl}(s),v_1)ds+   \int_{t_1}^{t}e^{-\nu_\epsilon(t-s)}(wg)(X_{cl}(s),v)ds,
    \end{aligned}
  \end{equation}
  for $v_{ 3}<0$, and
  \begin{equation}\label{2.2.6}
    h^{n+1}(x,v)=(1-\alpha)e^{-\nu_\epsilon(v)(t-t_1)}h^{n}(x_1,v_1)++e^{-\nu_\epsilon(v)(t-t_1)}\tilde{r}(v_0)+\int_{t_1}^{t}e^{-\nu_\epsilon(t_0-s)}(wg)(X_{cl}(s),v)ds,
  \end{equation}
  for $v_{ 3}>0$. Here we use the notation $\nu_\epsilon=\nu(v)+\epsilon.$
  By \eqref{2.2.5} and \eqref{2.2.6}, one can obtain that
  \begin{equation}\label{h-n}
  \left\|h^{n+1}\right\|_{L^\infty}\leq (1-\alpha)|h^n|_{L^\infty(\gamma)}+|w\tilde{r}|_{L^\infty(\gamma_-)}+\left\|\nu^{-1}wg\right\|_{L^\infty}.
  \end{equation}
  We now consider \eqref{2.2.4} with $n=0$. Since $f^0\equiv 0$, it is straightforward to see that
  \begin{equation}\label{2.2.8}
      \left\|h^1\right\|_{L^\infty}+|h^1|_{L^\infty(\gamma)} \leq C\{\left\|\nu^{-1}wg\right\|_{L^\infty}+|w \tilde{r}|_{L^\infty (\gamma_-)}\}.
  \end{equation}
  Therefore, we have solved \eqref{2.2.4} for $n=0$. Assume that we have already solved \eqref{2.2.4} for $0,1,2\dots,n$, and have obtained that
  \begin{equation}\label{2.2.9}
    \left\|h^{n}\right\|_{L^\infty}+|h^n|_{L^\infty(\gamma)} \leq C_{n}\{\left\|\nu^{-1}wg\right\|_{L^\infty }+|w \tilde{r}|_{L^\infty (\gamma_-)}\}.
  \end{equation}
  Then we can solve \eqref{2.2.4} with $n+1$ by using \eqref{2.2.5},\eqref{2.2.6}. We obtain the $L^\infty_{x,v}$ estimates for $n+1$ from \eqref{2.2.9} inductively:
  \begin{equation}\label{2.2.10}
  \begin{aligned}
      \left\|h^{n+1}\right\|_{L^\infty}+|h^{n+1}|_{L^\infty(\gamma)}
      &\leq (1-\alpha)|h^n|_{L^\infty(\gamma)}+C_n\{|w\tilde{r}|_{L^\infty(\gamma_-)}+\|\nu^{-1}wg||_{L^\infty}\}\\
      &\leq C_{n}\{\left\|\nu^{-1}wg\right\|_{L^\infty}+|w \tilde{r}|_{L^\infty (\gamma_-)}\} <\infty.
  \end{aligned}
  \end{equation}
  Therefore, we have solved \eqref{2.2.4} for all $n\in \mathbb{N}$ with \eqref{2.2.9}.
  
  Note that the constant $C_{n}$ depends on $n$, so we cannot take limit as $n\to \infty$. Similar to the derivation of \eqref{h-n}, we also obtain
  \begin{equation*}
    \begin{aligned}
      &\left\|h^{n+1}-h^n\right\|_{L^\infty}+|h^{n+1}-h^n|_{L^\infty(\gamma_+)}\leq (1-\alpha)(\left\|h^{n}-h^{n-1}\right\|_{L^\infty}+|h^{n}-h^{n-1}|_{L^\infty(\gamma_+)})\\
      &\leq \dots \leq (1-\alpha)^n(\left\|h^1\right\|_{L^\infty }+|h^1|_{L^\infty(\gamma_+)})\leq (1-\alpha)^n(\left\|\nu^{-1}wg\right\|_{L^\infty}+|w\tilde{r}|_{L^\infty}).
    \end{aligned}
  \end{equation*}
  Thus, we immediately conclude that the sequence $\{h^n\}^\infty_{n=0}$ is a Cauchy sequence in $L^\infty_{x,v}$, due to the accommodation coefficient $\alpha \in (0,1)$. Consequently, there exists a unique $h \in L^\infty_{x,v}$ such that
  $$  \left\|h^{n}-h\right\|_{L^\infty }+|h^n-h|_{L^\infty(\gamma_+)} \to 0 \quad   \text{as}\quad \! n \to \infty,$$
   with the bound
  \begin{equation*}
    \|h \|_{L^\infty} \leq C (\|\nu^{-1}wg\|_{L^\infty}+|w\tilde{r}|_{L^\infty}).
  \end{equation*}
  Therefore, we obtain a solution to \eqref{Con-eq1} with $a=0$:
  \begin{equation*}
    \begin{cases}
    \mathscr{L} _0 f =\epsilon f + v_3\partial_xf+\nu(v)f=g,\\
    f(0,v)|_{v_3>0}=(1-\alpha)L\gamma_+f+\tilde{r}(v),\\
  f(d,v)|_{v_3<0}= L \gamma _+ f(d, v).
    \end{cases}
    \end{equation*}
  with $L^2_{x,v}$ and $L^\infty_{x,v}$ weighted bound.
  
  \textbf{Step 2. Existence of \eqref{Con-eq1} in $0 < a \ll 1$.} 
  
   First we define the Banach space
  $$\begin{aligned}
    \mathbb{X}:=\{&f(x,v):wf \in L^\infty(\Omega_d\times \mathbb{R}^3),wf\in L^\infty(\gamma),\\
     &\text{ and } f(0,v)|_{v_3>0}=(1-\alpha)f(0,R_0 v)+\tilde{r}(v), f(d,v)|_{v_3<0}=f(d,R_d v)\}.
  \end{aligned}$$
  By the a priori estimates  Lemma \ref{Lmm-apriori-infty} and Lemma \ref{Lmm-apriori-L2}, one has that
  \begin{equation}\label{2.2.16}
    \left\|w \mathscr{L}^{-1}_ag\right\|_{L^\infty_{x,v}}+|w \mathscr{L}^{-1}_ag|_{L^\infty(\gamma)}=\left\|h^\epsilon\right\|_{L^\infty}+|h^\epsilon|_{L^\infty(\gamma)} \leq  \tfrac{C_{\epsilon,d}}{\alpha}\left(\left\|\nu^{-1}wg\right\|_{L^\infty}+|w\tilde{r}(v)|_{L^\infty(\gamma_-)}\right).
  \end{equation}
   Next, let $f^\epsilon_1= \mathscr{L}^{-1}_ag_1$ and $f^\epsilon_2=\mathscr{L}^{-1}_ag_2$, where $\nu^{-1}wg_1,\nu^{-1}wg_2\in L^{\infty}_{x,v}$ are two solutions of \eqref{Con-eq1} with $g$ replaced by $g_1$ and $g_2$, respectively. Then we have the following equation
   \begin{equation}\label{2.2.17}
    \begin{cases}
      \epsilon(f^\epsilon_2-f^\epsilon_1)+v_3\partial_x(f^\epsilon_2-f^\epsilon_1)+\nu(v)(f^\epsilon_2-f^\epsilon_1)-aK((f^\epsilon_2-f^\epsilon_1))=g_2-g_1,\\
      (f^\epsilon_2-f^\epsilon_1)(0,v)|_{v_3>0}=(1-\alpha)L\gamma_+(f^\epsilon_2-f^\epsilon_1)(0,v),\\
      (f^\epsilon_2-f^\epsilon_1)(0,v)|_{v_3<0}= L\gamma_+(f^\epsilon_2-f^\epsilon_1)(d,v).
    \end{cases}
   \end{equation}
  From Lemmas \ref{Lmm-apriori-infty} and \ref{Lmm-apriori-L2}, we infer that
  \begin{equation}\label{L2-difference}
    \left\|\mathscr{L}_a^{-1}g_1-\mathscr{L}_a^{-1}g_2\right\|^2_{L^2 }\leq \tfrac{C_\epsilon}{\alpha}\left\|g_1-g_2\right\|^2_{L^2},
  \end{equation}
  and
  \begin{equation}\label{Linfty-difference}
    \left\|w \mathscr{L}^{-1}_ag_1-w \mathscr{L}^{-1}_ag_2\right\|_{L^\infty }+|w \mathscr{L}^{-1}_ag_1-w \mathscr{L}^{-1}_ag_2|_{L^\infty(\gamma)}\leq \tfrac{C_{\epsilon,d}}{\alpha}\left\|\nu^{-1}w(g_1-g_2)\right\|_{L^\infty }.
  \end{equation}
  Note that $C_{\epsilon,d}$ here does not depend on $a$.
  
  Now, we define the operator
  $$\mathcal{J}_af=\mathscr{L}_0 ^{-1}(aKf+g).$$
  For any $f_1,f_2 \in \mathbb{X}$, using \eqref{Linfty-difference}, we obtain
  \begin{equation*}
    \begin{aligned}
      &\left\|w(\mathcal{J}_af_1-\mathcal{J}_af_2)\right\|_{L^\infty }+|w(\mathcal{J}_af_1-\mathcal{J}_af_2)|_{L^\infty(\gamma)}\\
      =&\left\|w\{\mathscr{L} ^{-1}_0(aKf_1+g)-\mathscr{L}^{-1}_0(aKf_2+g)\}\right\|_{L^\infty}\\
      &+|w\left(\mathscr{L} ^{-1}_0(aKf_1+g)-\mathscr{L} ^{-1}_0(aKf_2+g)\right)|_{L^\infty(\gamma)}\\
      \leq & a\tfrac{C_{\epsilon,d}}{\alpha}\left\|\nu^{-1}w(Kf_1-Kf_2)\right\|_{L^\infty }\\
      \leq & a\tfrac{C_{\epsilon,d}}{\alpha}\left\|w(f_1-f_2)\right\|_{L^\infty }.
    \end{aligned}
    \end{equation*}
  We choose $a_*$ small enough such that $a \tfrac{C_{\epsilon,d}}{\alpha}\leq \frac{1}{2}$, and hence, $\mathcal{J}_a:\mathbb{X} \to \mathbb{X} $ is a contraction mapping for $a \in [0,a_*]$. Therefore, by the Banach Fix-Point Theorem, there exists a unique fixed point $f\in\mathbb{X}$ such that $\mathcal{J}_af=f$, which immediately yields that
  $$\mathscr{L}_af=g \text{ for } a\in [0,a_*].$$
  Thus we have obtained the existence of \eqref{Con-eq1} for $a \in [0,a_*]$.
   
  \textbf{Step 3. Existension from $a=0$ to $a=1$.}
  
   Now, we define the operator
  \begin{equation*}
    \mathcal{J}_{a_*+a}=\mathscr{L}^{-1}_{a_*}(aKf+g).
  \end{equation*}
  Noting that the estimates for $ \mathscr{L} $ are independent of $a_*$. By a similar argument as in Step 2, we can prove that $\mathcal{J}_{a_*+a}:\mathbb{X}\to \mathbb{X}$ is also a contraction mapping for $a \in [0,a_*]$. Step by step, we get the existence result for the operator $\mathscr{L}^{-1}_a$ from $a=0$ to $a=1$. Finally, we observe that $\mathscr{L}_1^{-1}g$ satisfies the estimates \eqref{2.2.16}. The uniqueness follows from \eqref{L2-difference}. Therefore, we have completed the proof of Lemma \ref{Lmm-Con-eq}.
\end{proof}

\subsection{Uniform in $\epsilon$ estimates: limits from \eqref{Ap-eq-eps} to \eqref{Ap-eq-d}}
Lemma \ref{Lmm-Con-eq} establishes the existence of \eqref{Ap-eq-eps}. Now that we have a solution for \eqref{Ap-eq-eps}, to take the limit with respect to $\epsilon$, we require uniform estimates in $\epsilon$. The most challenging task is obtaining estimates for the macroscopic part. Motivated by \cite{esposito2013non,huang2022boundary,huang2023boundary}, we derive the uniform estimates in $\epsilon$ for the macroscopic part by selecting specific test functions in the weak formulation of \eqref{Ap-eq-eps}.
\begin{lemma}\label{Lmm-macro}
   let $d \geq 1$, and $f^\epsilon$ be the solution of \eqref{Ap-eq-eps} constructed in Lemma \ref{Lmm-Con-eq}. Then the following estimate holds
   \begin{equation}\label{macro}
    \left\|\mathbf{P}f^\epsilon\right\|^2_{L^2 }\leq C_d\{\left\|(\mathbf{I-P})f^\epsilon\right\|^2_\nu+\left\|g\right\|^2_{L^2 }+(2-\alpha)|f^\epsilon(0)|^2_{L^2\gamma_+}+|\tilde{r}|^2_{L^2(\gamma_-)}\},
   \end{equation}
   where $C_d>0$ is a constant independent of $\epsilon$.
\end{lemma}
\begin{proof}We denote
  \begin{equation*}
    \mathbf{P} f^\epsilon(x,v) =\left(a^\epsilon(x)+b^\epsilon(x) \cdot v+ c^\epsilon(x) (\tfrac{|v|^2-3}{2}) \right)\sqrt{m(v)}.
  \end{equation*}
   The weak formulation of \eqref{Ap-eq-eps} is
\begin{equation}\label{weak-form}
  \begin{aligned}
    &\epsilon\int_{0}^{d}\int_{\mathbb{R}^3}f^\epsilon(x,v)\psi dxdv-\int_{0}^{d}v_3f^\epsilon(x,v)\partial_x\psi(x,v)dxdv\\
  =&-\int_{\mathbb{R}^3}v_3f^\epsilon(d,v)\psi(d,v)dv+\int_{\mathbb{R}^3}v_3f^\epsilon(0,v)\psi(0,v)dv-\int_{0}^{d}\int_{\mathbb{R}^3}\psi(x,v)\mathcal{L}f^\epsilon(x,v) dx dv\\
  &+\int_{0}^{d}\int_{\mathbb{R}^3}g(x,v)\psi(x,v)dxdv.
  \end{aligned}
\end{equation}
As in \cite{huang2023boundary,huang2022boundary}, we choose a suitable test function $\psi$ to obtain the unifrom estimates in $\epsilon$ for the macroscopic part of $f^\epsilon$. 

We define
\begin{equation}\label{2.3.3}
  \zeta^\epsilon_c(x)=\int_{x}^{d}c^\epsilon(z)dz.
\end{equation}
It is straightforward to obtain
\begin{equation}\label{2.3.4}
  \zeta^\epsilon_c(d)=0, \quad |\zeta^\epsilon_c(0)|\leq \sqrt{d}\left\|c\right\|_{L^2_x},\quad \text{and}\quad \! \left\|\zeta^\epsilon_c(x)\right\|_{L^2_x}\leq d\left\|c^\epsilon\right\|_{L^2_x}.
\end{equation}
We define the test function $\psi$ in \eqref{weak-form} as
\begin{equation*}
  \psi=\psi_c=v_3(|v|^2-5)\sqrt{m}\zeta^\epsilon_c(x).
\end{equation*}
Then, the second term on left-hand side of \eqref{weak-form} is
\begin{equation}\label{2.3.6}
  \begin{aligned}
    &-\int_{0}^{d}\int_{\mathbb{R}^3}v_3f(x,v)\partial_x\psi(x,v)dxdv\\
=&\int_{0}^{d}\int_{\mathbb{R}^3}v_3^2\left(\tfrac{|v|^2-3}{2}\right)(|v|^2-5)c^\epsilon(x)^2m(v) dxdv\\
&+\int_{0}^{d} \int_{\mathbb{R}^3}(\mathbf{I-P})f^\epsilon(x,v)c^\epsilon(x)v_3^2(|v|^2-5)\sqrt{m}dxdv\\
\geq &5\int_{0}^{d}c^\epsilon(x)^2dx-C\left\|(\mathbf{I-P})f\right\|_{L^2_\nu}\left\|c\right\|_{L^2_x}\\
\geq &4\left\|c^\epsilon\right\|_{L^2_x}^2-C\left\|(\mathbf{I-P})f^\epsilon\right\|^2_{L^2_\nu},
  \end{aligned}
\end{equation}
where we have used the oddness and the fact that
\begin{equation*}
  \int_{\mathbb{R}^3}v_3^2(|v|^2-3)(|v|^2-5)m(v)dv=10,\quad  \int_{\mathbb{R}^3}v_3^2(|v|^2-5)m(v)dv=0.
\end{equation*}
Note that $\psi\in \mathcal{N}^\perp$. The first term on the left-hand side of \eqref{weak-form} can be estimated as
\begin{equation}\label{2.3.7}
\begin{aligned}
    \epsilon\int_{0}^{d}\int_{\mathbb{R}^3}f^\epsilon(x,v)\psi dxdv
    &=\epsilon\int_{0}^{d}\int_{\mathbb{R}^3}(\mathbf{I-P})f^\epsilon(x,v)v_3(|v|^2-5)\sqrt{m}\zeta^\epsilon_c(x) dxdv\\
  &\leq C\epsilon d\left\|(\mathbf{I-P})f^\epsilon\right\|_{\nu}\left\|c^\epsilon\right\|_{L^2_x}.
\end{aligned}
\end{equation}
For the boundary term on the right-hand side of \eqref{weak-form}, using \eqref{2.3.4} and the boundary condition in \eqref{Ap-eq-eps}, one has that
\begin{equation*}
  \begin{aligned}
    &-\int_{\mathbb{R}^3}v_3f^\epsilon(d,v)\psi(d,v)dv+\int_{\mathbb{R}^3}v_3f^\epsilon(0,v)\psi(0,v)dv\\
    = &\int_{v_3<0}v_3^2f^\epsilon(0,v)(|v|^2-5)\sqrt{m}\zeta^\epsilon_c(0)dv+\int_{v_3>0}v_3^2(|v|^2-5)\sqrt{m}[(1-\alpha)f(0,Rv)+\tilde{r}(v)]\zeta^\epsilon_c(0)dv\\
    \leq &C\sqrt{d}\left((2-\alpha)|f^\epsilon(0)|_{L^2(\gamma_+)}+|\tilde{r}(v)|_{L^2(\gamma_-)}\right)\left\|c\right\|_{L^2_x},
  \end{aligned}
\end{equation*}
 Therefore, the right-hand side of \eqref{weak-form} can be bounded by
\begin{equation}\label{2.3.9}
  d\{\left\|(\mathbf{I-P})f\right\|_{L^2_\nu}+\left\|g\right\|_{L^2}+(2-\alpha)|f^\epsilon(0)|_{L^2(\gamma_+)}+|\tilde{r}| _{L^2(\gamma_-)}\}\left\|c^\epsilon\right\|_{L^2_x}.
\end{equation}
By combining \eqref{2.3.6}, \eqref{2.3.7}, and \eqref{2.3.9} together, we finally obtain the estimate for $c^\epsilon$:
\begin{equation}\label{2.3.10}
  \left\|c^\epsilon\right\|^2_{L^2_x}\leq Cd^2\{\left\|(\mathbf{I-P})f^\epsilon\right\|^2_{L^2_\nu}+\left\|g\right\|^2_{L^2}+(2-\alpha)|f^\epsilon(0)|^2_{L^2(\gamma_+)}+|\tilde{r}|^2 _{L^2(\gamma_-)}\}.
\end{equation}

The corresponding estimates for $b^\epsilon(x)$ and $a^\epsilon(x)$ can be obtained by further consider 
$$\begin{aligned}
  &\zeta^\epsilon_{b,i}= \int_{x}^{d} b^\epsilon_i(z) dz,\quad  \psi_{b_i}=|v|^2v_3v_i \zeta^\epsilon_{b,i}(x),\quad \! (i=1,2),\\
  &\zeta^\epsilon_{b,3}= \int_{x}^{d} b^\epsilon_3(z) dz,\quad  \psi_{b_i}=(v_3^2-1)\sqrt{m} \zeta^\epsilon_{b,3}(x),
\end{aligned}$$
and
$$\zeta^\epsilon_{a}= \int_{x}^{d} a^\epsilon(z) dz,\psi_{a}=v_3(|v|^2-10)\sqrt{m} \zeta^\epsilon_{a}(x).$$
The proof is similar, so we omit the details for simplicity. Therefore, we have completed the proof of Lemma \ref{Lmm-macro}.
\end{proof}
Based on Lemma \ref{Lmm-macro}, we can obtain the uniform estimates in $\epsilon$ and take the limit as $\epsilon\to 0$.
\begin{lemma}\label{lm2.4}
   Let $d\geq 1$, under the assumption in Theorem \ref{Thm-linear}, there exists a unique solution $f=f(x,v)$ to the linearized steady Boltzmann equation
   \begin{equation}\label{Ap-eq-d1}
    \begin{cases}
v_3\partial_xf+Lf=g,  \quad (x,v)\in \Omega_d \times \mathbb{R}^3,\\
    f(0,v)|_{v_3>0}=(1-\alpha)L\gamma_+f(0,v)+\tilde{r}(v),\\
  f(d,v)|_{v_3<0}=L\gamma_+f(d,v),
    \end{cases}
    \end{equation}
    with
    \begin{equation}\label{2.4.4}
      \begin{aligned}
        &\left\|f \right\|_{L^2}^2+\alpha|f |^2_{L^2(\gamma_+)}\leq \tfrac{C_{d}}{\alpha}(\left\|g\right\|^2_{L^2}+|\tilde{r}(v)|^2_{L^2(\gamma_-)}),\\
        &\left\|wf\right\|_{L^\infty }+|wf(0)|_{L^\infty(\gamma_-)} \leq \tfrac{C_d}{\alpha}(\left\|\nu^{-1}w g\right\|_{L^\infty }+|w\tilde{r}|_{L^\infty_{v}(\gamma_-)}).
      \end{aligned}
    \end{equation}
\end{lemma}
\begin{proof}
   Multiplying \eqref{Ap-eq-eps} by $f^\epsilon$ and integrating respect to $x$ and $v$, one obtains that
  \begin{equation*}
    2\epsilon\left\|f^\epsilon\right\|_{L^2}^2+\alpha|f^\epsilon|^2_{L^2(\gamma_+)}+2c_0\left\|(\mathbf{I-P})f^\epsilon\right\|^2_{L^2_\nu}\leq 2\left\|f^\epsilon\right\|_{L^2}
  \left\|g\right\|_{L^2}+ \tfrac{C}{\alpha}|\tilde{r}(v)|^2_{\gamma_-},
  \end{equation*}
  which together with Lemma \ref{Lmm-macro}, we have
  \begin{equation}\label{2.4.2}
  \left\|f^\epsilon\right\|_{L^2}^2+\alpha|f^\epsilon|^2_{L^2(\gamma_+)}+2c_0\left\|(\mathbf{I-P})f^\epsilon\right\|^2_{L^2_\nu}\leq \tfrac{C_{d}}{\alpha}(\left\|g\right\|^2_{L^2}+|\tilde{r}(v)|^2_{L^2(\gamma_-)}).
  \end{equation}
  By Lemma \ref{Lmm-apriori-infty} and \eqref{2.4.2}, we get that
  \begin{equation}\label{2.4.5}
    \left\|wf^\epsilon\right\|_{L^\infty}+|wf^\epsilon|_{L^\infty(\gamma)}\leq \tfrac{C_{d}}{\alpha} (
    \left\|\nu^{-1}wg\right\|_{L^\infty}+|w\tilde{r}(v)|_{L^\infty(\gamma_-)}).
  \end{equation}
  Now, we consider the limit of $f^\epsilon$ as $\epsilon \to 0$. Consider $f^{\epsilon_1},f^{\epsilon_2}$, the difference $f^{\epsilon_1}-f^{\epsilon_2}$ satisfies the following system:
  \begin{equation*}
    \begin{cases}
      v_3\partial_x(f^{\epsilon_1}-f^{\epsilon_2})+\mathcal{L}(f^{\epsilon_1}-f^{\epsilon_2})=-\epsilon_1f^{\epsilon_1}+\epsilon_2f^{\epsilon_2},\\
      (f^{\epsilon_1}-f^{\epsilon_2})(0,v)|_{v_3>0}=(1-\alpha)L\gamma_+(f^{\epsilon_1}-f^{\epsilon_2}),\\
      (f^{\epsilon_1}-f^{\epsilon_2})(d,v)|_{v_3<0}=L\gamma_+(f^{\epsilon_1}-f^{\epsilon_2}).
    \end{cases}
  \end{equation*}
  By similar arguments above, one has
  \begin{equation*}
    \begin{aligned}
      &\left\|(f^{\epsilon_1}-f^{\epsilon_2})\right\|_{L^2}^2+\left\|(\mathbf{I-P})(f^{\epsilon_1}-f^{\epsilon_2})\right\|_\nu^2+\alpha|(f^{\epsilon_1}-f^{\epsilon_2})|_{L^2(\gamma_+)}^2\\
    &\leq \tfrac{C_{d}}{\alpha}\left\|\epsilon_1f^{\epsilon_1}-\epsilon_2f^{\epsilon_2}\right\|^2_{L^2}\\
    &\leq \tfrac{C_{d}}{\alpha}(\epsilon_1+\epsilon_2)^2 (\left\|\nu^{-1}wg\right\|^2_{L^\infty}+|w\tilde{r}|^2_{L^\infty(\gamma_-)}),
    \end{aligned}
  \end{equation*}
  and
  \begin{equation*}
    \left\|w(f^{\epsilon_2}-f^{\epsilon_1})\right\|_{L^\infty}+|w(f^{\epsilon_1}-f^{\epsilon_2})|_{L^\infty(\gamma)} \leq \tfrac{C_{d}}{\alpha}(\epsilon_1+\epsilon_2)\left(\left\|\nu^{-1}wg\right\|_{L^\infty}+|wr|_{L^\infty(\gamma_-)}\right).
  \end{equation*}
  Therefore, the sequence $f^\epsilon$ is a Cauchy sequence in $L^\infty_{x,v}$ as $\epsilon \to 0$, and there exists a solution $f$ to  \eqref{Ap-eq-d1}  satisfies \eqref{2.4.4}. The uniqueness follows from \eqref{2.4.4}. Thus, we have completed the proof of Lemma \ref{lm2.4}.
\end{proof}
We claim that the system \eqref{Ap-eq-d1} is equivalent to \eqref{Ap-eq-d}. To verify this, it suffices to show that 
\begin{equation*}
  P_\gamma f(0,v)= \int_{v_3>0}v_3 f\sqrt{m}dv =\lambda.
\end{equation*}
Multiplying  \eqref{Ap-eq-d1} by $\sqrt{m}$ and integrating respect over $(x,v)$ and noting that $g\in \mathcal{N}^\perp$, we obtain
\begin{equation*}
  \int_{\R^3}v_3 f(d,v) \sqrt{m} dv -\int_{\R^3}v_3 f(0,v) dv=0.
\end{equation*}
The first term on the left-hand side vanishes due to the boundary condition at $x=d$. Thus, we have
\begin{multline*}
  \int_{\R^3}v_3 f(0,v) \sqrt{m} dv =\int_{v_3<0}v_3 f(0,v) \sqrt{m} dv \\+\int_{v_3>0}v_3 ((1-\alpha)f(0,R_0v)+ \alpha \sqrt{2\pi m(v)} \lambda +r(v)) \sqrt{m} dv=0.
\end{multline*}
Hence,
\begin{equation*}
  P_\gamma f(0,v)= \int_{v_3>0}v_3 f\sqrt{m}dv =\lambda.
\end{equation*}
Therefore, we conclude that the system \eqref{Ap-eq-d1} is equivalent to \eqref{Ap-eq-d}.
\subsection{Uniform in $d$ estimates}\label{Sec-uniformd}
In this subsection, we aim to obtain the uniform estimates in $d$ for the equation \eqref{Ap-eq-d}, which is equivalent to \eqref{Ap-eq-d1}, so that we can take the limit $d \to \infty$.

Let $f(x,v)$ be the solution of \eqref{Ap-eq-d1}, and define the macroscopic projection
\begin{equation*}
  \mathbf{P}f(x,v)=\left(a(x)+b(x)\cdot v+c(x)(\tfrac{|v|^2-3}{2})\right)\sqrt{m}.
\end{equation*}
We multiply \eqref{Ap-eq-d1} by $\sqrt{m}$ and integrate respect to $v$ over $\mathbb{R}^3$, yielding
\begin{equation*}
\frac{d}{dx}  \int_{\mathbb{R}^3}v_3\sqrt{m}f(x,v)dv=\int_{\mathbb{R}^3}g(x,v)\sqrt{m}dv\equiv 0,
\end{equation*}
since $g \in \mathcal{N}^\perp$. From the specular reflection boundary condition at $x=d$, we have
\begin{equation*}
  \int_{\mathbb{R}^3}v_3f(d,v)\sqrt{m}dv=0.
\end{equation*}
Therefore, it follows that
\begin{equation}\label{2.5.3}
  b_3(x) \equiv 0  \quad \text{for all}\quad \! x\in [0,d].
\end{equation}
This results in
\begin{equation}\label{2.5.4}
  \int_{\mathbb{R}^3}v_3|\mathbf{P}f(x,v)|^2dv \equiv0.
\end{equation}
Similarly, multiplying \eqref{Ap-eq-d1} by $v_i\sqrt{m}, i=1,2,$ and $\tfrac{|v|^2-5}{2}\sqrt{m}$, we can obtain the following relations
\begin{equation}\label{AB=0}
  \begin{aligned}
    &\int_{\mathbb{R}^3}A_{i3}   (\mathbf{I-P})f=0 \quad \text{ for } i=1,2,\\
  &\int_{\mathbb{R}^3}B_3 (\mathbf{I-P})f=0,
  \end{aligned}
\end{equation}
where $A_{ij},B_i$ are introduced in \eqref{AB}. Using $b_3(x)\equiv 0$, the above expressions allow us to conclude
\begin{equation}\label{2.5.6}
  \int_{\mathbb{R}^3}v_3\mathbf{P}f(x,v) (\mathbf{I-P})f(x,v)\equiv 0, \quad \text{for all}\quad\! x\in [0,d].
\end{equation}
From \eqref{2.5.4} and \eqref{2.5.6}, we conclude that
\begin{equation}\label{2.5.7}
  \int_{\mathbb{R}^3}v_3|f(x,v)|^2dv= \int_{\mathbb{R}^3}v_3|(\mathbf{I-P})f(x,v)|^2dv ,\quad \forall x\in [0,d].
\end{equation}
These relations will be highly useful in deriving the uniform estimates in $d$ for the macroscopic part, which is an essential step in proving the desired result as $d\to \infty$.

To obtain exponential decay estimates for the solution $f(x,v)$, we introduce a new function $\tilde{f}(x,v)$ as follows
\begin{equation*}
  \tilde{f}(x,v)=f(x,v)-q(v),
\end{equation*}
where $q(v)$ is a macroscopic function, given by
\begin{equation*}
  q(v)=\left(q_0+q_1v_1+q_2v_2+q_3 (\tfrac{|v|^2-3}{2}) \right)\sqrt{m},
\end{equation*}
with the constants $q_0,q_1,q_2,q_3$ to be determined. The system governing $\tilde{f}$ can be written as:
\begin{equation}\label{tilde-f}
  \begin{cases}
v_3\partial_x\tilde{f}+L\tilde{f}=g,  \\
  \tilde{f}(0,v)|_{v_3>0}=(1-\alpha)L\gamma_+\tilde{f}(0,v)+\tilde{r}(v)-\alpha q(v),\\
\tilde{f}(d,v)|_{v_3<0}=\tilde{f}(d,Rv).
  \end{cases}
  \end{equation}
Moreover, there holds 
\begin{equation*}
  (\mathbf{I-P})\tilde{f}(x,v)=(\mathbf{I-P})f(x,v),
\end{equation*}
and
\begin{equation*}
    \mathbf{P}\tilde{f}=\left(\tilde{a}+\tilde{b}_1(x)v_1+\tilde{b}_2(x)v_2+\tilde{c}(x)(\tfrac{|v|^2-3}{2})\right)\sqrt{m(v)}, 
\end{equation*}
with
\begin{equation}\label{2.6.5}
  \begin{cases}
    \tilde{a}(x)=a(x)-q_0,\\
    \tilde{b}_i(x)=b_i ( x )-q_i, \quad i=1,2, \\
    \tilde{c}(x)=c(x)-q_3.
  \end{cases}
\end{equation} 

As in \cite{huang2023boundary,huang2022boundary}, we choose $q$ as follows:
\begin{lemma}\label{lm2.6}
  There exist constants $q_{0}, q_{1}, q_{2}, q_{3}$ such that
\begin{equation}\label{q0-q3}
  \begin{aligned}
  & \int_{\mathbb{R}^{3}} v_{3} \tilde{f}(d, v) \cdot v_{3} \sqrt{m} d v=0, \\
  & \int_{\mathbb{R}^{3}} v_{3} \tilde{f}(d, v) \cdot \mathcal{L}^{-1}\left( A _{3 i}\right) d v=0, \quad i=1,2, \\
  & \int_{\mathbb{R}^{3}} v_{3} \tilde{f}(d, v) \cdot \mathcal{L}^{-1}\left( B _{3}\right) d v=0,
  \end{aligned}
\end{equation}
where the constants $q_{0}, q_{1}, q_{2}, q_{3}$ are determined by $\left(  {a}_{1},   {b}_{1},   {b}_{2},   {c}\right)(d)$ and $(\mathbf{I}-\mathbf{P})  f (d)$.
\end{lemma}
\begin{proof}
  As in Lemma 2.11 in \cite{huang2022boundary}, \eqref{2.6.5} and \eqref{q0-q3} imply that
  \begin{equation}\label{define-q}
    \left(\begin{array}{cccc}
    1 & 0 & 0 & 1 \\
    0 & \kappa_{1} & 0 & 0 \\
    0 & 0 & \kappa_{1} & 0 \\
    0 & 0 & 0 & \kappa_{2}
    \end{array}\right)\left(\begin{array}{c}
    q_{0} \\
    q_{1} \\
    q_{2} \\
    q_{3}
    \end{array}\right)=\left(\begin{array}{c}
      {a}(d)+  {c}(d)+\int_{\mathbb{R}^{3}}(\mathbf{I}-\mathbf{P})   {f}(d, v) \cdot  A _{33}(v) d v \\
    \kappa_{1}   {b}_{1}(d)+\int_{\mathbb{R}^{3}} v_{3}(\mathbf{I}-\mathbf{P})   {f}(d, v) \cdot \mathcal{L}^{-1}\left( A _{31}\right) d v \\
    \kappa_{1}   {b}_{2}(d)+\int_{\mathbb{R}^{3}} v_{3}(\mathbf{I}-\mathbf{P})   {f}(d, v) \cdot \mathcal{L}^{-1}\left( A _{32}\right) d v \\
    \kappa_{2}   {c}(d)+\int_{\mathbb{R}^{3}} v_{3}(\mathbf{I}-\mathbf{P})   {f}(d, v) \cdot \mathcal{L}^{-1}\left( B _{3}\right) d v
    \end{array}\right) ,
  \end{equation}
  where $\kappa_1$ and $\kappa_2$ are positive constants, as defined in \eqref{kappa}. Since the matrix is nonsingular, the values of $\left(q_{0}, q_{1}, q_{2}, q_{3}\right)$ are determined. Thus, the proof of Lemma \ref{lm2.6} is complete.
\end{proof}

The following lemma provides the $L^2_{x,v}$ estimates for $e^{\sigma x}\tilde{f}$. 
\begin{lemma}\label{Lmm-tilde-f}
  Let $q_{0}, q_{1}, q_{2}, q_{3}$ be the constants determined in Lemma \ref{lm2.6}. Then it holds that
\begin{equation}\label{2.7.1}
  \left\|e^{\sigma x} \tilde{f}\right\|_{L ^{2}} \leq \tfrac{C}{\alpha} \left(|\tilde{r}|_{L^{2}\left(\gamma_{-}\right)}+\tfrac{1}{\sigma_{1}-\sigma}\left\|e^{\sigma_{1} x} g\right\|_{L ^{2}}\right),
\end{equation}
with $0<\sigma<\sigma_{1} \leq \sigma_{0}$, and $C>0$ is a constant depending only on $\sigma_{1}$ and $\alpha$.
\end{lemma}

\begin{proof}First, we estimate the microscopic part. Note that $(\mathbf{I-P})\tilde{f}=(\mathbf{I-P})f$, so we directly estimate $(\mathbf{I-P})f$. We claim that for $\sigma_1\in(0,\sigma_0]$, it holds that
\begin{equation}\label{claim-1}
  \alpha |f(0)|^2_{L^2(\gamma_+)}+\int_{0}^{d}e^{2\sigma_1x}\left\|(\mathbf{I-P})f(x)\right\|_\nu^2dx   \leq  \tfrac{C}{\alpha} ( |\tilde{r}|^2_{L^2(\gamma_-)}+ \|e^{\sigma_1 x}g \|^2_{L^2}) ,
\end{equation}
Indeed, multiplying \eqref{Ap-eq-d1} by $e^{2\sigma_1x}f$ and integrating over $(x,v)\in \Omega_d \times \R^3$, one obtains that
\begin{multline}\label{A1}
  \frac{d}{dx}\iint_{\Omega_d \times\mathbb{R}^3}v_3e^{2\sigma_1x}f^2dxdv-\iint_{\Omega_d \times \mathbb{R}^3}v_3 2\sigma_1e^{2\sigma_1x}f^2dxdv\\
  + \iint_{\Omega_d\times \mathbb{R}^3}e^{2\sigma_1x}f  \L f dx dv 
  =\iint_{\Omega_d \times \mathbb{R}^3}e^{2\sigma_1x}fg dxdv.
\end{multline}
The boundary conditions in $\eqref{Ap-eq-d1}$ infer that
\begin{equation}\label{A2}
  \begin{aligned}
    &\frac{d}{dx}\iint_{\Omega_d \times\mathbb{R}^3}v_3e^{2\sigma_1x}f^2dxdv =-\int_{\R^3}v_3 f^2(0,v) dv\\
    =& -\int_{v_3>0} \big(v_3(1-\alpha)f(0,Rv)+\tilde{r}(v) \big)^2dv-\int_{v_3<0}v_3f^2(0,v)dv\\
     \geq &-\big(1-(1-\alpha)^2\big)\int_{v_3<0}v_3 f^2(0,v)dv - 2(1-\alpha)\int_{v_3>0}v_3 f(0, R_0 v) \tilde{r}(v) dv \\
     \geq &-\alpha\int_{v_3<0}v_3 f^2 (0,v) d v-\tfrac{C}{\alpha} \int_{v_3>0}v_3\tilde{r}^2(v)dv.
  \end{aligned}
\end{equation}
Using the relation \eqref{2.5.7} and the coercivity of $\L$, one obtains
\begin{equation}\label{A3}
  -\iint_{\Omega_d \times \mathbb{R}^3}v_3 2\sigma_1e^{2\sigma_1x}f^2dxdv
  + \iint_{\Omega_d\times \mathbb{R}^3}e^{2\sigma_1x}f  \L f dx dv \geq (c_0-C\sigma_1)\|e^{\sigma_1 x}(\mathbf{I-P})f\|^2_\nu,
\end{equation}
for some $c_0>0$. Using the fact that $g\in \mathcal{N}^\perp$, one obtains 
\begin{equation}\label{A4}
  \iint_{\Omega_d \times \mathbb{R}^3}e^{2\sigma_1x}fgdv \leq \eta \|e^{\sigma_1 x}(\mathbf{I-P})f\|^2_\nu +C_\eta \|e^{\sigma_1 x}g\|^2_{L^2}.
\end{equation}
Using \eqref{A1}-\eqref{A4} and taking $\sigma_0,\eta $ suitably small, one can conclude the claim \eqref{claim-1}.

Next, we consider the macroscopic part. Multiplying $\eqref{tilde-f}$ by $\mathcal{L}^{-1}\left( A _{31}\right), \mathcal{L}^{-1}\left( A _{32}\right)$, and $\mathcal{L}^{-1}\left( B _{3}\right)$, respectively, we obtain
\begin{equation}\label{2.7.2}
  \begin{aligned}
  & \frac{d}{d x}\left(\begin{array}{c}
  \int_{\mathbb{R}^{3}} v_{3} \tilde{f}(x, v) \mathcal{L}^{-1}\left( A _{31}\right) d v \\
  \int_{\mathbb{R}^{3}} v_{3} \tilde{f}(x, v)  \mathcal{L}^{-1}\left( A _{32}\right) d v \\
  \int_{\mathbb{R}^{3}} v_{3} \tilde{f}(x, v)  \mathcal{L}^{-1}\left( A _{3}\right) d v
  \end{array}\right)=\left(\begin{array}{c}
  \int_{\mathbb{R}^{3}}g  \mathcal{L}^{-1}\left( A _{31}\right) d v \\
  \int_{\mathbb{R}^{3}}g  \mathcal{L}^{-1}\left( A _{32}\right) d v \\
  \int_{\mathbb{R}^{3}}g  \mathcal{L}^{-1}\left( B _{3}\right) d v
  \end{array}\right).
  \end{aligned}
\end{equation}
Here we use the fact that $(\mathbf{I-P})\tilde{f}=(\mathbf{I-P})f$, together with the relation \eqref{AB=0}. Integrating the above system over $[x, d]$ and using the relation \eqref{q0-q3}, we obtain
\begin{equation*}
  \left(\begin{array}{c}
  \int_{\mathbb{R}^{3}} v_{3} \tilde{f}(x, v)  \mathcal{L}^{-1}\left( A _{31}\right) d v \\
  \int_{\mathbb{R}^{3}} v_{3} \tilde{f}(x, v)  \mathcal{L}^{-1}\left( A _{32}\right) d v \\
  \int_{\mathbb{R}^{3}} v_{3} \tilde{f}(x, v)  \mathcal{L}^{-1}\left( B _{3}\right) d v
  \end{array}\right)=-\int_{x}^{d}\left(\begin{array}{c}
  \int_{\mathbb{R}^{3}} g \mathcal{L}^{-1}\left( A _{31}\right) d v \\
  \int_{\mathbb{R}^{3}} g \mathcal{L}^{-1}\left( A _{32}\right) d v \\
  \int_{\mathbb{R}^{3}} g \mathcal{L}^{-1}\left( B _{3}\right) d v
  \end{array}\right)(z) d z.
\end{equation*}
Then we have
\begin{equation}\label{2.7.5}
  \left|\left(\kappa_{1} \tilde{b}_{1}, \kappa_{1} \tilde{b}_{2}, \kappa_{2} \tilde{c}\right)(x)\right| \leq C\|(\mathbf{I}-\mathbf{P})   {f}(x)\|_{\nu}+C \int_{x}^{d}\|g(z)\|_{L_{v}^{2}} d z.
\end{equation}
Multiplying \eqref{2.7.5} by $e^{\sigma x}$ with $0<\sigma<\sigma_{1} \leq \sigma_{0}$ and using \eqref{claim-1}, one has
\begin{equation}\label{2.7.6}
  \begin{aligned}
  & \int_{0}^{d} e^{2 \sigma x}\left|\left(\tilde{b}_{1}, \tilde{b}_{2}, \tilde{c}\right)(x)\right|^{2} d x \\
   \leq & C \int_{0}^{d} e^{2 \sigma x}\|(\mathbf{I}-\mathbf{P})   {f}(x)\|_{\nu}^{2} d x+C \int_{0}^{d} e^{2 \sigma x}\left|\int_{x}^{d}\|g(z)\|_{L_{v}^{2}} d z\right|^{2} d x \\
  \leq & \tfrac{C}{\alpha} \left(|\tilde{r}|_{L^{2}\left(\gamma_{-}\right)}^{2}+\tfrac{1}{\sigma_{1}\left(\sigma_{1}-\sigma\right)}\left\|e^{\sigma_{1} x} g\right\|_{L ^{2}}^{2}\right).
  \end{aligned}
\end{equation}
Finally, we consider the estimate for $\tilde{a}$. Specifically, multiplying \eqref{tilde-f} by $v_{3} \sqrt{m}$ gives
\begin{equation*}
  \frac{d}{d x} \int_{\mathbb{R}^{3}} \tilde{f}(x, v) v_{3}^{2} \sqrt{m} d v=\int_{\mathbb{R}^{3}} v_3g   \sqrt{m} d v=0.
\end{equation*}
Integrating above equation over $[x, d]$, and using $\eqref{q0-q3}_1$, one obtains
\begin{equation}\label{2.7.8}
  \tilde{a}(x)=-\tilde{c}(x)+\int_{\mathbb{R}^{3}}(\mathbf{I}-\mathbf{P})   {f}(x, v) v_3^2 \sqrt{m} d v.
\end{equation}
Multiplying \eqref{2.7.8} by $e^{\sigma x}$ with $0<\sigma<\sigma_{1} \leq \sigma_{0}$ and using \eqref{claim-1} and \eqref{2.7.6}, it holds that
\begin{equation}\label{2.7.9}
  \int_{0}^{d} e^{2 \sigma x}|\tilde{a}(x)|^{2} d x \leq \tfrac{C}{\alpha}\left(|\tilde{r}|_{L^{2}\left(\gamma_{-}\right)}^{2}+\tfrac{1}{\sigma_{1}\left(\sigma_{1}-\sigma\right)}\left\|e^{\sigma_{1} x} g\right\|_{L ^{2}}^{2}\right).
\end{equation}
We have proved \eqref{2.7.1} by using \eqref{claim-1}, \eqref{2.7.6} and \eqref{2.7.9}. Therefore, Lemma \ref{Lmm-tilde-f} is established.
\end{proof}

\begin{corollary}
  Let $\tilde{f}$ be the solution of \eqref{tilde-f}, then for $0<\sigma<\sigma_0$, it holds that
\begin{equation}\label{2.7.10}
  \left\|e^{\sigma x} w \tilde{f}\right\|_{L^{\infty}}+\left|e^{\sigma x} w \tilde{f}\right|_{L^{\infty}(\gamma)} 
   \leq \tfrac{C}{\alpha }\left( \tfrac{1}{\sigma_{0}-\sigma}\left\|e^{\sigma_{0} x} \nu^{-1} w g\right\|_{L ^{\infty}}+\left|\left(q_{0}, q_{1}, q_{2}, q_{3}\right)\right|+|w \tilde{r}|_{L^{\infty}\left(\gamma_{-}\right)}\right),
\end{equation}
for some constant $C>0$.
\end{corollary}
 \begin{proof}
  Let $\tilde{h}:=e^{\sigma x} w \tilde{f}$. Multiply \eqref{tilde-f} by $e^{\sigma x} w$ to obtain
 \begin{equation*}
   \begin{cases}
     v_{3} \partial_{x} \tilde{h}+  (\nu(v)-\sigma v_3) \tilde{h}= K_{w} \tilde{h}+e^{\sigma x} w g, \\
      \tilde{h}(0, v)|_{v_{3}>0}=(1-\alpha)L\gamma_+\tilde{h}(0,v)-\alpha w(v)q+\tilde{r}(v),\\
     \tilde{h}(d, v)|_{v_{3}<0}= L\gamma_+\tilde{h}(d, v).
   \end{cases}
 \end{equation*}
 We take $\sigma_{0}>0$ small such that $\nu(v)-\sigma v_3 \geq \frac{1}{2} \nu(v)>0$. By the same arguments as in Lemma \ref{Lmm-apriori-infty} and using Lemma \ref{Lmm-tilde-f}, we can obtain
 $$
 \begin{aligned}
 & \|\tilde{h}\|_{L^{\infty}}+|\tilde{h}|_{L^{\infty}(\gamma)} \\
 \leq & C\left(\left\|e^{\sigma x} \tilde{f}\right\|_{L^{2}}+\left\|e^{\sigma x} \nu^{-1} w g\right\|_{L^{\infty}}+\left|\left(q_{0}, q_{1}, q_{2}, q_{3}\right)\right|+|w \tilde{r}|_{L^{\infty}\left(\gamma_{-}\right)}\right) \\
 \leq & \tfrac{C}{\alpha}\left(\tfrac{C}{\sigma_{1}-\sigma}\left\|e^{\sigma_{1} x} g\right\|_{L^{2}}+\left\|e^{\sigma x} \nu^{-1} w g\right\|_{L ^{\infty}}+\left|\left(q_{0}, q_{1}, q_{2}, q_{3}\right)\right|+|w\tilde{r}|_{L^{\infty}\left(\gamma_{-}\right)}\right) \\
 \leq & \tfrac{C}{\alpha}\left(\tfrac{1}{\sigma_{0}-\sigma}\left\|e^{\sigma_{0} x} \nu^{-1} w g\right\|_{L ^{\infty}}+\left|\left(q_{0}, q_{1}, q_{2}, q_{3}\right)\right|+|w\tilde{r}|_{L^{\infty}\left(\gamma_{-}\right)}\right),
 \end{aligned}
 $$
 where we have chosen $\sigma_{1}=\frac{\sigma_{0}+\sigma}{2}$ such that $0<\sigma<\sigma_{1}<\sigma_{0}$.
 \end{proof}
Note that the estimate in \eqref{2.7.10} depends on $q_i$. To take limit as $d \to \infty$, we need to study the asymptotic behavior of $q_i(d)(i=0,1,2,3)$. Specifically, we first prove the uniform boundness of $q^\infty$ and obtain the uniform estimates for $\tilde{f}$.
\begin{lemma}\label{lm2.9}
 Let $\tilde{f}$ be the solution of \eqref{tilde-f}. Then, we have the following.
  \begin{equation}\label{2.9.1}
    |(q_0,q_1,q_2,q_3)(d)| \leq C_\alpha\left(\left\|e^{\sigma_0x}\nu^{-1}wg\right\|_{L^\infty}+|w\tilde{r}|_{L^\infty(\gamma_-)}\right),
  \end{equation}
  and
  \begin{equation}\label{2.9.9}
    \begin{aligned}
    & \left\|e^{\sigma x} w \tilde{f}\right\|_{L ^{\infty}}+\left|e^{\sigma x} w \tilde{f}\right|_{L^{\infty}(\gamma)} \leq C_\alpha\left(\tfrac{1}{\sigma_{0}-\sigma}\left\|e^{\sigma_{0} x} \nu^{-1} w g\right\|_{L ^{\infty}}+|w \tilde{r}|_{L^{\infty}\left(\gamma_{-}\right)}\right),
    \end{aligned}
  \end{equation}
   for some $C_\alpha>0$.
\end{lemma}
\begin{proof}
   Recall \eqref{define-q} the definition of $q_i(i=0,\cdots,3)$. We first estimate the macroscopic quantities $a,b_1,b_2,c$.
  Multiplying \eqref{Ap-eq-d1} by $\mathcal{L}^{-1}(A_{13})$ and integrating over $[0,x]$, one obtains that
  \begin{equation*}
    \int_{\mathbb{R}^3}v_3f(x,v)\mathcal{L}^{-1}(A_{13})dv=\int_{\mathbb{R}^3}v_3\mathcal{L}^{-1}(A_{13})f(0,v)dv+\int_{0}^{x}\int_{\mathbb{R}^3}g(z)\mathcal{L}^{-1}A_{13}dzdv,
  \end{equation*}
  where we have used the fact that $\int_{\R^3} A_{13}(\mathbf{I-P})f dv=0$. A direct calculation shows that:
  \begin{equation}\label{2.9.2}
  \begin{aligned}
      \kappa_1b_1(x)=&-\int_{\mathbb{R}^3}v_3\mathcal{L}^{-1}(A_{13})(\mathbf{I-P})f(x,v)dv+\int_{\mathbb{R}^3}v_3\mathcal{L}^{-1}(A_{13})f(0,v)dv\\
      &+\int_{0}^{x}\int_{\mathbb{R}^3}g(z)\mathcal{L}^{-1}A_{13}dzdv.
  \end{aligned}
  \end{equation}
  Using the boundary condition $\eqref{Ap-eq-d1}_2$ and \eqref{2.4.4}, the boundary term on the right-hand side of \eqref{2.9.2} can be bounded by
  \begin{equation*}
    \int_{\mathbb{R}^3}v_3\mathcal{L}^{-1}(A_{13})f(0,v)dv\leq C\left(|f(0,v)|_{L^2(\gamma_+)}+|\tilde{r}|_{L^2(\gamma_-)}\right)
    \leq C_\alpha( \|g \|_{L^2}+| \tilde{r}|_{L^2(\gamma_-)}).
  \end{equation*}
  Therefore, one obtains that
  \begin{equation}\label{2.9.4}
    |b_1(x)| \leq C \|(\mathbf{I-P})f(x,v)\|_{L^\infty_v}+C_\alpha(\left\|g\right\|_{L^2}+| \tilde{r}|_{L^2(\gamma_-)}),
  \end{equation}
  where we have used \eqref{claim-1}. Similarly, multiplying \eqref{Ap-eq-d1} by $\mathcal{L}^{-1}A_{23},\mathcal{L}^{-1}B_{3}$ and $v_3\sqrt{m}$, respectively, we can obtain
  \begin{equation}\label{2.9.10}
    \begin{aligned}
        \kappa_1b_2(x)=&-\int_{\mathbb{R}^3}v_3\mathcal{L}^{-1}A_{23}(\mathbf{I-P})f(x,v)dv+\int_{\mathbb{R}^3}v_3\mathcal{L}^{-1}A_{23}f(0,v)dv\\
        &+\int_{0}^{x}\int_{\mathbb{R}^3}g(z)\mathcal{L}^{-1}A_{23}dzdv.
    \end{aligned}
    \end{equation}
      \begin{equation}\label{2.9.11}
        \begin{aligned}
            \kappa_2c(x)=&-\int_{\mathbb{R}^3}v_3\mathcal{L}^{-1}B_{3}(\mathbf{I-P})f(x,v)dv+\int_{\mathbb{R}^3}v_3\mathcal{L}^{-1} B_{3}f(0,v)dv\\
            &+\int_{0}^{x}\int_{\mathbb{R}^3}g(z)\mathcal{L}^{-1}B_{3}dzdv.
        \end{aligned}
        \end{equation}
        \begin{equation}\label{2.9.12}
          \begin{aligned}
              a(x)=&-c(x)-\int_{\mathbb{R}^3}v_3 A_{33}(\mathbf{I-P})f(x,v)dv+\int_{\mathbb{R}^3}v_3^2 f(0,v)\sqrt{m}dv,
          \end{aligned}
          \end{equation}
  and the estimates for $a,b_2$ and $c$:
  \begin{equation}\label{2.9.5}
    |(a,b_2,c)(x)| \leq C \|(\mathbf{I-P})f(x,v)\|_{L^\infty_v}+C_\alpha(\left\|g\right\|_{L^2}+| \tilde{r}|_{L^2(\gamma_-)}).
  \end{equation}
  Note that the constants $C_\alpha$ in \eqref{2.9.4} and \eqref{2.9.5} are independent of $d$. Then from \eqref{2.9.4}, \eqref{2.9.5} and \eqref{2.7.10}, we can bound $q_i$ by
  \begin{equation}\label{2.9.6}
  \begin{aligned}
   \left|\left(q_{0}, q_{1}, q_{2}, q_{3}\right)(d)\right|
    &\leq C\left|\left(  {a},   {b}_{ 1},   {b}_{2},   {c}\right)(d)\right|+C\|(\mathbf{I}-\mathbf{P})   {f}(d)\|_{L_{v}^{\infty}} \\
  & \leq C \|(\mathbf{I}-\mathbf{P}) \tilde{f}(d) \|_{L_{v}^{\infty}}+C_\alpha\left(\|g\|_{L ^{2}}+| \tilde{r}|_{L^{2}\left(\gamma_{-}\right)}\right) \\
  & \leq C e^{-\sigma d}\|e^{\sigma d} w \tilde{f}(d)\|_{L_{v}^{\infty}}+C_\alpha\left(\|g\|_{L ^{2}}+| \tilde{r}|_{L^{2}\left(\gamma_{-}\right)}\right) \\
  & \leq C_\alpha e^{-\sigma d}\left|\left(q_{0}, q_{1}, q_{2}, q_{3}\right)(d)\right|+C_\alpha\left(\tfrac{1}{\sigma_{0}-\sigma}\left\|e^{\sigma_{0} x} \nu^{-1} w g\right\|_{L ^{\infty}}+|w  \tilde{r}|_{L^{\infty}\left(\gamma_{-}\right)}\right).
  \end{aligned}
  \end{equation}
  Let $\sigma=\frac{1}{2} \sigma_{0}$, and choose $d_{0}$ sufficiently large such that $C_\alpha e^{-\sigma d_{0}} \leq \frac{1}{2}$. Then \eqref{2.9.6} reduces to
  \begin{equation}\label{2.9.7}
    \left|\left(q_{0}, q_{1}, q_{2}, q_{3}\right)(d)\right| \leq C_\alpha\left(\left\|e^{\sigma_{0} x} \nu^{-1} w g\right\|_{L ^{\infty}}+|w  \tilde{r}|_{L^{\infty}\left(\gamma_{-}\right)}\right) \quad \text {for} \quad \! d \geq d_{0} \geq 1.
  \end{equation}
  For $d \in\left[1, d_{0}\right]$, using \eqref{define-q} and \eqref{2.4.4}, we have
  \begin{equation*} 
    \left|\left(q_{0}, q_{1}, q_{2}, q_{3}\right)(d)\right| \leq C_\alpha\left(\left\|\nu^{-1} w g\right\|_{L ^{\infty}}+|w  \tilde{r}|_{L^{\infty}\left(\gamma_{-}\right)}\right),
  \end{equation*}
  which, together with \eqref{2.9.7}, leads to \eqref{2.9.1}. Finally, \eqref{2.9.9} follows directly from \eqref{2.7.10} and \eqref{2.9.1}. Thus, the proof of Lemma \ref{lm2.9} is complete.
\end{proof}
In the following, we study the asymptotic behavior of $\left(q_{0}, q_{1}, q_{2}, q_{3}\right)(d)$ as $d \to \infty$. Let $f_{d_1}$ and $f_{d_2}$ be solutions of \eqref{Ap-eq-d1} with $1 \leq d_{1} \leq d_{2}<\infty$. Then it holds that
\begin{equation}\label{2.10.1}
  \left\{\begin{array}{l}
  v_{3} \partial_{x}\left(  {f}_{d_{2}}-  {f}_{d_{1}}\right)+  \mathcal{L}\left({f}_{d_{2}}-{f}_{d_{1}}\right)=0, \quad(x, v) \in\left(0, d_{1}\right) \times \mathbb{R}^{3}, \\
  \left({f}_{d_{2}}-{f}_{d_{1}}\right)(0, v)|_{v_{3}>0}=(1-\alpha)L\gamma_+\left(  {f}_{d_{2}}-  {f}_{d_{1}}\right)(0, v).
  \end{array}\right.
\end{equation}
\begin{lemma}\label{Lmm-apriori-infty0}
 Let $d \geq1$. It holds that
\begin{equation}\label{2.10.0}
   \left|\left(  {f}_{d_{2}}-  {f}_{d_{1}}\right)(0)\right|_{L^{2}\left(\gamma_{+}\right)} \leq C_\alpha e^{-\sigma_{1} d_{1}}\left(| \tilde{r}|_{L^{2}\left(\gamma_{-}\right)}^{2}+\int_{0}^{d_{1}} e^{2 \sigma_{1} x}\|g(x)\|_{L_v ^{2}}^{2} d x \right)^{\frac{1}{2}},
\end{equation}
for some $C_\alpha> 0$.
 \end{lemma}
\begin{proof}
  Multiplying \eqref{2.10.1} by $  {f}_{d_{2}}-  {f}_{d_{1}}$ and integrating over $[0,x]$  with $x\in (0,d_1]$, we obtain
  \begin{multline*}
    - \frac{1}{2} \int_{\mathbb{R}^{3}} v_{3}\left|\left(  {f}_{d_{2}}-  {f}_{d_{1}}\right)(0, v)\right|^{2} d v +c_{0}   \int_{0}^{x}\left\|(\mathbf{I}-\mathbf{P})\left(  {f}_{d_{2}}-  {f}_{d_{1}}\right)(z)\right\|_{\nu}^{2} d z  \\
     \leq -\frac{1}{2} \int_{\mathbb{R}^{3}} v_{3}\left|\left(  {f}_{d_{2}}-  {f}_{d_{1}}\right)(x, v)\right|^{2} d v,
  \end{multline*}
 where we have used the coercivity of $\L$. Next, using the boundary condition in $\eqref{2.10.1}$, we have
  \begin{equation*}
    -\int_{\mathbb{R}^{3}} v_{3}\left|\left(  {f}_{d_{2}}-  {f}_{d_{1}}\right)(0, v)\right|^{2} d v=(2-\alpha)\alpha|(f_{d_2}-f_{d_1})|^2_{L^2(\gamma_+)}.
  \end{equation*}
  Moreover, the relations \eqref{2.5.3} and \eqref{AB=0} infer that
  \begin{equation*}
    \begin{aligned}
    & \int_{\mathbb{R}^{3}} v_{3}\left|\left(  {f}_{d_{2}}-  {f}_{d_{1}}\right)(x, v)\right|^{2} d v \\
     =&\int_{\mathbb{R}^{3}} v_{3}\left|  {f}_{d_{1}}(x, v)\right|^{2} d v-2 \int_{\mathbb{R}^{3}} v_{3}   {f}_{d_{2}}(x, v)   {f}_{d_{1}}(x, v) d v+\int_{\mathbb{R}^{3}} v_{3}\left|  {f}_{d_{2}}(x, v)\right|^{2} d v \\
     =&\int_{\mathbb{R}^{3}} v_{3}\left|(\mathbf{I}-\mathbf{P})   {f}_{d_{1}}(x, v)\right|^{2} d v-2 \int_{\mathbb{R}^{3}} v_{3}(\mathbf{I}-\mathbf{P})   {f}_{d_{2}}\left(d_{1}, v\right)  (\mathbf{I}-\mathbf{P})   {f}_{d_{1}}\left(d_{1}, v\right) d v \\
    & +\int_{\mathbb{R}^{3}} v_{3}\left|(\mathbf{I}-\mathbf{P})   {f}_{d_{2}}(x, v)\right|^{2} d v.
    \end{aligned}
  \end{equation*}
  As a result, one obtains
  \begin{equation}\label{2.10.5}
    |{f}_{d_{2}}-  {f}_{d_{1}}(0)|_{L^{2}\left(\gamma_{+}\right)}^{2} \leq C_{\alpha}\left(\left\|(\mathbf{I}-\mathbf{P})   {f}_{d_{1}}(x)\right\|_{\nu}^{2}+\left\|(\mathbf{I}-\mathbf{P})   {f}_{d_{2}}(x)\right\|_{\nu}^{2}\right), \quad \forall x \in\left[0, d_{1}\right].
  \end{equation}
  Integrating \eqref{2.10.5} over $x \in\left[d_{1}-1, d_{1}\right]$, and using \eqref{claim-1}, we have
  \begin{equation}\label{2.10.6}
    \begin{aligned}
    & \left|\left(  {f}_{d_{2}}-  {f}_{d_{1}}\right)(0)\right|_{L^{2}\left(\gamma_{+}\right)}^{2} \\
     \leq &C_\alpha\left(\int_{d_{1}-1}^{d_{1}}\left\|(\mathbf{I}-\mathbf{P})   {f}_{d_{1}}(x)\right\|_{\nu}^{2} d x+\int_{d_{1}-1}^{d_{1}}\left\|(\mathbf{I}-\mathbf{P})   {f}_{d_{2}}(x)\right\|_{\nu}^{2} d x\right) \\
     \leq &C_\alpha e^{-2 \sigma_{1} d_{1}}\left(\int_{d_{1}-1}^{d_{1}} e^{2 \sigma_{1} x}\left\|(\mathbf{I}-\mathbf{P})   {f}_{d_{1}}(x)\right\|_{\nu}^{2} d x+\int_{d_{1}-1}^{d_{1}} e^{2 \sigma_{1} x}\left\|(\mathbf{I}-\mathbf{P})   {f}_{d_{2}}(x)\right\|_{\nu}^{2} d x\right) \\
     \leq &C_\alpha e^{-2 \sigma_{1} d_{1}}\left(|\tilde{r}|_{L^{2}\left(\gamma_{-}\right)}^{2}+\int_{0}^{d_{1}} e^{2 \sigma_{1} x}\|g(x)\|_{L_v ^{2}}^{2} d x\right),
    \end{aligned}
  \end{equation}
  which immediately yields \eqref{2.10.0}. Therefore the proof of Lemma \ref{Lmm-apriori-infty0} is complete.
\end{proof}

\begin{lemma}\label{Lmm-q-asymptotic}
  There exist constants $\left(q_{0}^{\infty}, q_{1}^{\infty}, q_{2}^{\infty}, q_{3}^{\infty}\right)$ such that
\begin{equation}\label{2.11.1}
  \lim _{d \to+\infty}\left(q_{0}, q_{1}, q_{2}, q_{3}\right)(d)=\left(q_{0}^{\infty}, q_{1}^{\infty}, q_{2}^{\infty}, q_{3}^{\infty}\right),
\end{equation}
with
\begin{equation}\label{2.11.2}
  \left|\left(q_{0}^{\infty}, q_{1}^{\infty}, q_{2}^{\infty}, q_{3}^{\infty}\right)\right| \leq C\left(\left\|e^{\sigma_{0} x} \nu^{-1} w g\right\|_{L ^{\infty}}+|w \tilde{r}|_{L^{\infty}\left(\gamma_{-}\right)}\right).
\end{equation}
\end{lemma}
\begin{proof}
   Letting $1 \leq d_{1} \leq d_{2}<\infty$, it follows from \eqref{define-q}, \eqref{2.7.10}, \eqref{2.9.2}, \eqref{2.9.10}-\eqref{2.9.12} and \eqref{2.10.0},  that
  \begin{equation}\label{2.11.3}
    \begin{aligned}
    & \left|\left(q_{0}\left(d_{2}\right)-q_{0}\left(d_{1}\right), q_{1}\left(d_{2}\right)-q_{1}\left(d_{1}\right), q_{2}\left(d_{2}\right)-q_{2}\left(d_{1}\right), q_{3}\left(d_{2}\right)-q_{3}\left(d_{1}\right)\right)\right| \\
    \leq &\left|\left(  {a}_{d_{2}}\left(d_{2}\right)-  {a}_{d_{1}}\left(d_{1}\right),   {b}_{d_{2}, 1}\left(d_{2}\right)-  {b}_{d_{1}, 1}\left(d_{1}\right),   {b}_{d_{2}, 2}\left(d_{2}\right)-  {b}_{d_{1}, 2}\left(d_{1}\right),   {c}_{d_{2}}\left(d_{2}\right)-  {c}_{d_{1}}\left(d_{1}\right)\right)\right| \\
    &  +C\|(\mathbf{I}-\mathbf{P}) \tilde{f}_{d_{1}}\left(d_{1}\right)\|_{L_{v}^{\infty}}+C\|(\mathbf{I}-\mathbf{P}) \tilde{f}_{d_{2}}\left(d_{2}\right)\|_{L_{v}^{\infty}} \\
     \leq &\left|\left(  {f}_{d_{2}}-  {f}_{d_{1}}\right)(0)\right|_{L^{2}\left(\gamma_{+}\right)} +\int_{d_{1}}^{d_{2}}\|g(z)\|_{L_{v}^{2}} d z   +C_\alpha e^{-\sigma d_{1}}\left(\tfrac{1}{\sigma_{0}-\sigma}\left\|e^{\sigma_{0} x} \nu^{-1} w g\right\|_{L ^{\infty}}+|w \tilde{r}|_{L^{\infty}\left(\gamma_{-}\right)}\right) \\
    \leq  &C_\alpha e^{-\sigma d_{1}}\left(\tfrac{1}{\sigma_{0}-\sigma}\left\|e^{\sigma_{0} x} \nu^{-1} w g\right\|_{L ^{\infty}}+|w \tilde{r}|_{L^{\infty}\left(\gamma_{-}\right)}\right) \to 0, \quad \text {as} \quad\! d_{1} \to+\infty,
    \end{aligned}
  \end{equation}
  which immediately gives \eqref{2.11.1}. Furthermore, \eqref{2.11.2} follows directly from Lemma \ref{lm2.9}. Therefore the proof of Lemma \ref{Lmm-q-asymptotic} is complete.
\end{proof}

\subsection{Limits from \eqref{Ap-eq-d} to \eqref{Auxi-eq}: Proof of Lemma \ref{Lmm-auxi}} Let $(\tilde{f}_{d_1}, q (d_1))$ and $(\tilde{f}_{d_2},q_(d_2))$ be the solutions of \eqref{tilde-f} with $1 \leq d_{1} \leq d_{2}<\infty$. We denote
$$
\bar{f}(x, v):=\left(\tilde{f}_{d_{2}}-\tilde{f}_{d_{1}}\right)(x, v) \quad \text { and} \quad \!\bar{h}(x, v)=w(v)  \bar{f}(x, v), \quad \text{for} \quad \!(x, v) \in\left[0, d_{1}\right] \times \mathbb{R}^{3}.
$$
Then it follows from \eqref{tilde-f} that
\begin{equation}\label{2.12.1}
  \left\{\begin{array}{l}
  v_{3} \partial_{x}  \bar{f}+  \mathcal{L}  \bar{f}=0, \quad(x, v) \in \Omega_{d_1} \times \mathbb{R}^{3}, \\
   \bar{f}(0, v)|_{v_{3}>0}=(1-\alpha) L\gamma_+ \bar{f}(0, v)-\alpha\left(q(d_2)-q(d_1)\right).
  \end{array}\right.
\end{equation}
Using the estimates from subsection \ref{Sec-uniformd} and following similar arguements as those in subsection 3.3 of \cite{huang2023boundary}, we derive the following $L^2$ and $L^\infty$ estimates:
\begin{equation}\label{2.12.6}
  |\bar{f}(x,v)|^2_{L^2(\gamma_+)}+ \iint_{\Omega_{d_1} \times\mathbb{R}^{3}}| \bar{f}(x, v)|^{2} d v d x \leq C_\alpha d_{1} e^{-2 \sigma d_{1}}\left(\tfrac{1}{\sigma_{0}-\sigma}\left\|e^{\sigma_{0} x} \nu^{-1} w g\right\|_{L ^{\infty}}+|w \tilde{r}|_{L^{\infty}\left(\gamma_{-}\right)}\right)^{2},
\end{equation}
and
\begin{equation}
   \|\bar{h}\|_{L^{\infty}\left(\left[0, \frac{1}{2} d_{1}\right] \times \mathbb{R}^{3}\right)}+|\bar{h}(0)|_{L^{\infty}\left(\gamma_{+}\right)} 
  \leq  C_\alpha d_1e^{-\frac{1}{2}\sigma d_1}  \left(\tfrac{1}{\sigma_0-\sigma}\left\|e^{\sigma_0 x}\nu^{-1}wg\right\|_{L^\infty}+|w\tilde{r}|_{L^\infty(\gamma_-)}\right).
\end{equation}
The proof is similar so we omit the details here for brevity. 
Hence, there exists a function $f(x,v)$ so that $ \left\|w(\tilde{f_d}-f)\right\|_{L^\infty}+\left|w(f(0)-\bar{f}(0))\right|_{L^\infty(\gamma)} \to 0$ as $d_1 \to \infty$. Based on the strong convergence and \eqref{2.9.9}, we have
\begin{equation}\label{2.12.20}
  \left\|w e^{\sigma x}f\right\|_{L^\infty}+\left| wf(0)\right|_{L^\infty(\gamma)}\leq C_\alpha \big( \tfrac{1}{\sigma_0-\sigma}\left\|e^{\sigma_0x}\nu^{-1}wg\right\|_{L^\infty}+\left|w\tilde{r}\right|_{L^\infty(\gamma_-)} \big) .
\end{equation}
It is direct to see $f$ satisfies
\begin{equation}\label{2.12.14}
\begin{cases}
    v_3\partial_x f+\mathcal{L}f=g, \quad (x,v)\in \mathbb{R}_+\times \mathbb{R}^3,\\
    f(0,v)|_{v_3>0}=(1-\alpha)L\gamma_+f-\alpha q^\infty(v)+\tilde{r}(v),\\
\end{cases}
\end{equation}
with
\begin{equation}\label{2.12.13}
  q^{\infty}(v)=\left(a^{\infty}+b_{1}^{\infty} v_{1}+b_{2}^{\infty} v_{2}+c^{\infty}(\tfrac{|v|^{2}-3}{2})\right) \sqrt{m(v)},
\end{equation}
  such that
\begin{equation}\label{2.12.21}
  \left|\left(a^{\infty}, b_{1}^{\infty}, b_{2}^{\infty}, c^{\infty}\right)\right| \leq C\left(\left\|e^{\sigma_{0} x} \nu^{-1} w g\right\|_{L ^{\infty}}+|w \tilde{r}(v)|_{L^{\infty}\left(\gamma_{-}\right)}\right).
  \end{equation}
Let $\hat{f}=f+q^\infty$. We claim that $P_\gamma \hat{f}=\lambda$. In fact, multiplying \eqref{2.12.14} by $\sqrt{m}$ and integrating over $(x,v)$, we have
\begin{equation*}
  \int_{\R^3}v_3 q^\infty dv-\int_{\R^3} v_3 \hat{f}(0,v) dv=0,
\end{equation*}
where we have used the fact that $g \in \mathcal{N}^\perp$. The oddness of $q^\infty$ implies that $\int_{\R^3}v_3 q^\infty dv=0$. Recall \eqref{tilde-r} the definition of $\tilde{r}$. Then,
\begin{equation*}
  \begin{aligned}
    0=\int_{\R^3} \hat{f}(0,v) dv= \int_{v_3<0}v_3 \hat{f}(0,v)dv +\int_{v_3>0} v_3\big( (1-\alpha)L\gamma_+ \hat{f}(0,v) +\tilde{r}(v)\big)  dv\\
    =\alpha \int_{v_3<0} v_3 \hat{f}(0,v) dv  +\alpha \int_{v_3>0} \lambda v_3\sqrt{2\pi} m(v) dv
    =-\alpha P_\gamma \hat{f}(0,v) +\alpha \lambda 
  \end{aligned}
\end{equation*}
where we have used the fact that $\int_{v_3>0} v_3 \sqrt{m(v)}r(v) dv=0$. Therefore, $\hat{f}$ is a solution to \eqref{Auxi-eq}. Moreover,
the estimates \eqref{1.26} and \eqref{1.27} directly follow from \eqref{2.12.20} and \eqref{2.12.21}, since $|w\tilde{r}|_{L^\infty_v(\gamma_+)}\leq |w\tilde{r}|_{L^\infty_v(\gamma_+)}+|\lambda|$.

To prove the uniqueness, we let $f_1$ and $f_2$ be two solutions of \eqref{Auxi-eq} satisfying \eqref{1.26} and \eqref{1.27}, then it holds that
\begin{equation}\label{2.13.1}
  \begin{cases}
    v_3\partial_x (f_1-f_2)+\mathcal{L}{(f_1-f_2)}=0,\quad (x,v)\in \mathbb{R}_+\times \mathbb{R}^3,\\
    (L^R-\alpha L^D)(f_1-f_2)=0,\\
    P_\gamma(f_1-f_2)(0,v)=0.
  \end{cases}
\end{equation}
Multiplying $\eqref{2.13.1}_1$ by $f_1-f_2$, and integrating over $\R_+\times \mathbb{R}^3$, using the boundary condition, we have
\begin{equation}
  -\frac{(2-\alpha)\alpha}{2}|f_1-f_2|^2_{L^2(\gamma_+)}+c_0\left\|(\mathbf{I-P})(f_1-f_2)\right\|^2_{\nu}\leq 0.
\end{equation}
This implies
\begin{equation*}
  (\mathbf{I-P})(f_1-f_2)\equiv 0 \quad \text{ and\qquad} (f_1-f_2)|_{\gamma} \equiv 0.
\end{equation*}
Multiplying $\eqref{2.13.1}_1$ by $\mathcal{L}^{-1}A_{13},\mathcal{L}^{-1}A_{23},\mathcal{L}^{-1}B_{3} $ and $v_3\sqrt{m}$, respectively, by the same arguments as in Lemma \ref{Lmm-tilde-f} one has that
\begin{equation*}
  \left\|\mathbf{P}(f_1-f_2)\right\|_{L^2}\equiv 0.
\end{equation*}
Consequently, we have $f_1 \equiv f_2.$ Therefore we have proved the uniqueness.
$\hfill\square$

\section{Proof of the main result}
\subsection{Existence of the linear problem \eqref{KL-eq}: Proof of Theorem \ref{Thm-linear}}\label{Sec-linear-KL}
Lemma \ref{Lmm-auxi} ensures that, for any given $\lambda\in \R$, there exists a unique solution $f_\lambda$ to \eqref{Auxi-eq} satisfying $\lim_{x\to \infty} f_\lambda = q^\infty_\lambda$, where
\begin{equation*}
 q^\infty_{\lambda}=\left(a^\infty_\lambda +b^\infty_{\lambda,1}v_1+b^\infty_{\lambda,2}v_2 +c^\infty_{\lambda}(\tfrac{|v|^2-3}{2})\right)\sqrt{m(v)}.
\end{equation*}
Thus, the pair $(\tilde{f}_\lambda, \tilde{q}^\infty_\lambda) :=(f_\lambda-q^\infty_\lambda,q^\infty_\lambda-a^\infty_\lambda \sqrt{m(v)})$ solves \eqref{KL-eq} with the following estimates:
\begin{equation}\label{3.1.20}
 \begin{array}{l}
   \left\|e^{\sigma x}w\tilde{f}_\lambda\right\|_{L^\infty}+| w\tilde{f}_\lambda|_{L^\infty(\gamma)} \leq C\left(\tfrac{1}{\sigma_0-\sigma}\left\|e^{\sigma_0x}\nu^{-1}wg\right\|_{L^\infty}+|w r(v)|_{L^\infty(\gamma_-)}+|\lambda|\right),\\
   \|w\tilde{q}^\infty_\lambda\|_{L^\infty_v}\leq C\left(\|e^{\sigma_0x}\nu^{-1}wg\|_{L^\infty}+|wr|_{L^\infty(\gamma_-)}+|\lambda|\right).
 \end{array}
\end{equation}
Note that these estimates depend on $\lambda$. To complete the proof of Theorem \ref{Thm-linear}, it remains to show that $\tilde{f}_\lambda$ and $\tilde{q}^\infty_\lambda$ are independent of $\lambda$, and to establish uniform estimates with respect to $\lambda$.

First, let $(f_0,q_0^\infty)$ be the solution to \eqref{Auxi-eq} with $\lambda=0$. Then, $(\tilde{f}_0, \tilde{q}^\infty_0)$ solves \eqref{KL-eq}
with estimates independent of $\lambda$:
\begin{equation}\label{3.1.1}
\begin{array}{l}
    \left\|e^{\sigma x}w\tilde{f}_0\right\|_{L^\infty}+| w\tilde{f}_0(0)|_{L^\infty(\gamma)} \leq C \left(\frac{1}{\sigma_0-\sigma}\left\|e^{\sigma_0x}\nu^{-1}wg\right\|_{L^\infty}+|w r(v)|_{L^\infty(\gamma_-)}\right) ,\\
    \|w \tilde{q}^\infty_0\|_{L^\infty_v}\leq C\left(\|e^{\sigma_0x}\nu^{-1}wg\|_{L^\infty}+|wr|_{L^\infty(\gamma_-)} \right).
\end{array}
\end{equation}
To prove independence of $(\tilde{f}_\lambda, \tilde{q}^\infty_\lambda)$ with respect to $\lambda$, let $(f_{\lambda_1},q^\infty_{\lambda_1})$ and $(f_{\lambda_2},q^\infty_{\lambda_2})$ be two solutions of \eqref{Auxi-eq} with $\lambda_1\neq\lambda_2$. Define
  $$\hat{f}=f_{\lambda_1}-f_{\lambda_2}.$$
The system for $\hat{f}$ is
\begin{equation}\label{3.1.4}
  \begin{cases}
    v_3\partial_x \hat{f}+\mathcal{L}\hat{f}=0, \\
    (L^R-\alpha L^D)\hat{f}=0,\\
    P_\gamma \hat{f} = \lambda_1-\lambda_2,\\
    \lim_{x \to \infty}e^{\sigma x} \left\|w(v)[\hat{f}-q^\infty_{\lambda_1}+q^\infty_{\lambda_2}]\right\|_{L^\infty_v}=0.
  \end{cases}
\end{equation}
Multiplying $\eqref{3.1.4}_1$ by $\hat{f}$ and integrating over $(x,v)\in \mathbb{R}_+ \times \mathbb{R}^3$,
\begin{equation}\label{3.1.5}
\frac{1}{2}\int_{\mathbb{R}^3 }v_3|q_{\lambda_1}^\infty-q_{\lambda_2}^\infty|^2d v +c_0\int_{0}^{\infty } \left\|(\mathbf{I}-\mathbf{P})\hat{f}\right\|_\nu^2dx  \le
\frac{1}{2}\int_{\mathbb{R}^3 }v_3|\hat{f}(0, v )|^2dv.
\end{equation}
The oddness implies that
\begin{equation}\label{3.1.6}
\int_{\mathbb{R}^3 } v _3|q_{\lambda_1}^\infty-q_{\lambda_2}^\infty|^2d v =0.
\end{equation}
Using the Maxwell reflection boundary condition $\eqref{3.1.4}_2$, we obtain
\begin{equation}\label{3.1.7}
\begin{aligned}
\int_{\mathbb{R}^3 } v _3|\hat{f} (0, v )|^2d v =&\int_{ v _3>0 } v _3|\hat{f} (0, v )|^2d v +\int_{ v _3<0 } v _3|\hat{f} (0, v )|^2d v \\
=&(1-\alpha)^2\int_{ v _3>0 } v _3|\hat{f} (0,R_0 v )|^2d v + \sqrt{2\pi}\alpha^2 \left(\int_{ u_3<0 }|u_3|\hat{f} (0, u)\sqrt{m}d u\right) ^2\\
&+2\sqrt{2\pi}(1-\alpha)\alpha \int_{ v _3>0} v _3\hat{f} (0,R_0 v )\sqrt{m}d v \int_{ u_3<0 } |u_3|\hat{f} (0, u )\sqrt{m}d u\\
&+\int_{ v _3<0 } v _3|\hat{f} (0, v )|^2d v \\
=&(2-\alpha)\alpha \int_{ v _3<0 } v _3|\hat{f} (0, v )|^2d v +(2-\alpha)\alpha \sqrt{2\pi} \left(\int_{ v _3<0 } | v _3 |\hat{f} (0, v )\sqrt{m}d v \right)^2\\
=&-(2\alpha-\alpha^2)|(I-K\gamma_+)\hat{f} (0)|^2_{L^2{(\gamma_+)}}.
\end{aligned}
\end{equation}
 Substituting \eqref{3.1.6} and \eqref{3.1.7} into \eqref{3.1.5}, we obtain
\begin{equation}\label{3.1.8}
(2\alpha-\alpha^2)|(I-K\gamma_+)\hat{f} (0)|^2_{L^2{(\gamma_+)}}+\int_{0}^{\infty } \left\|(\mathbf{I}-\mathbf{P})\hat{f} \right\|^2_\upsilon dy \le 0,
\end{equation}
which implies
\begin{equation*} 
(\mathbf{I}-\mathbf{P})\hat{f} \equiv 0,
\end{equation*}
and
\begin{equation*} 
\hat{f} (0, v )=K\gamma_+ \hat{f} (0, v ) = (\lambda_1-\lambda_2)\sqrt{m}, \quad \text{for} \quad \!v_3<0.
\end{equation*}
Using the boundary condition $\eqref{3.1.4}_2$, we deduce
\begin{equation}
  \gamma_-\hat{f}=(1-\alpha)L\gamma_+\hat{f}+\alpha K\gamma_+\hat{f}=(\lambda_1-\lambda_2)\sqrt{m}.
\end{equation}
Thus,
\begin{equation}\label{3.1.11}
 \hat{f} (0, v )\equiv (\lambda_1-\lambda_2)\sqrt{m}.
\end{equation}
Next, we calculate $\mathbf{P} \hat{f}$. We denote
\begin{equation}\label{3.1.12}
\hat{f}(x) =\left(\hat{a}(x)+\hat{b}_1(x)v_1+\hat{b}_2(x)v_2+\hat{c}(x)(\tfrac{| v |^2-3}{2})\right)\sqrt{m}.
\end{equation}
  Multiplying $\eqref{3.1.4}_1$ by $ \mathcal{L}^{-1}B_3 $, and integrating over $ \mathbb{R}^3$ with respect to $ v $, using the fact that $\mathcal{L}$ is a self-adjoint operator, we obtain
  \begin{equation}\label{3.1.13}
    \frac{d}{dx}\int_{\mathbb{R}^3 }\hat{f} (x, v ) v _3\mathcal{L}^{-1}B_3d v +\int_{\mathbb{R}^3}\hat{f}(x,v) B_3( v )d v =0. \\
    \end{equation}
 Using \eqref{3.1.12}, we have $\int_{\mathbb{R}^3 } \hat {f}(x ,v)B_3( v ) d v =0,$ for $x \in [0,\infty)$. Therefore, the second term on the left-hand side of \eqref{3.1.13} vanishes. Consequently, for $x\in [0,\infty)$, we have
  \begin{equation*}
  \int_{\mathbb{R}^3 }\hat{f} (x, v ) v _3 \L^{-1}B_3d v =\int_{\mathbb{R}^3 }\hat{f} (0, v ) v _3 \L^{-1}B_3d v. \\
  \end{equation*}
  From \eqref{3.1.11}, one obtains that
  \begin{equation}\label{3.1.16}
    \int_{\mathbb{R}^3 }\hat{f} (0, v ) v _3 \L^{-1}B_3d v =0. \\
    \end{equation}
Using \eqref{3.1.12} and \eqref{3.1.16}, for any $x \in [0,\infty)$, one has
\begin{equation*}
  \int_{\mathbb{R}^3 }\hat{f} (x, v ) v _3 \L^{-1}B_3d v \equiv \int_{\mathbb{R}^3 }\hat{f} (0, v ) v _3 \L^{-1}B_3d v =0.
\end{equation*}
and therefore
\begin{equation}\label{3.1.18}
  \kappa_2 \hat{c}(x) \equiv 0.
\end{equation}
Similarly, multiplying $\eqref{3.1.4}_1$ by $\mathcal{L}^{-1}A_{13}$ and $\mathcal{L}^{-1}A_{23}$ and preforming the same arguments as above we can obtain that
\begin{equation*}
  \hat{b}_1(x)=\hat{b}_2(x) \equiv0.
\end{equation*}
As for $\hat{a}(x)$, we multiply $\eqref{3.1.4}_1$  by $v_3\sqrt{m}$ and obtain that:
\begin{equation*}
  \hat{a}(x)+\hat{c}(x)=P_\gamma\hat{f}(0,v).
\end{equation*}
Together with \eqref{3.1.18}, we obtain
\begin{equation*}
  \hat{a}(x)=P_\gamma\hat{f}(0,v).
\end{equation*}
Therefore, we have proved the fact that $\hat{f}(x,v)\equiv P_\gamma\hat{f}(0,v)\sqrt{m(v)}=(\lambda_1-\lambda_2)\sqrt{m}$. Together with $\eqref{3.1.4}_4$, we have
\begin{equation}\label{2.13.5}
  (q^\infty_{\lambda_1}-q^\infty_{\lambda_2})(v)=(\lambda_1-\lambda_2)\sqrt{m},
\end{equation}
and
\begin{equation*}
  f_{\lambda_1}-q^\infty_{\lambda_1}= f_{\lambda_2}-q^\infty_{\lambda_2}.
\end{equation*}
Therefore, we have proved that $\tilde{f}_\lambda \equiv\tilde{f}_0$ for any $\lambda\in \mathbb{R}$. Furthermore, \eqref{2.13.5} implies that $\tilde{q}^\infty_{\lambda} \equiv \tilde{q}^\infty_{0}$. Consequently, the solution $(\tilde{f},\tilde{q}^\infty)$ to  problem \eqref{1.1.20} is unique under the constraints $\eqref{1.1.19}_1$. Moreover, \eqref{1.1.19} and \eqref{1.1.20} follow from \eqref{3.1.1}. Therefore the proof of Theorem \ref{Thm-linear} is complete.

\subsection{Nonlinear problem: Proof of Theorem \ref{Thm-nonlinear}}
To solve the nonlinear problem \eqref{nonlinear-KL-eq}, we propose the following iteration system.
\begin{equation}\label{3.2.1}
\begin{cases}
    v_3\partial_x f_{i+1}+\mathcal{L}f_{i+1}=\Gamma(f_i,f_i)+S(x,v),\\
(L^R-\alpha L^D)f_{i+1}=-(L^R-\alpha L^D)q^\infty_{i+1}+R(v),\\
\lim_{x\to \infty}f_{i+1}=0,
\end{cases}
\end{equation}
where $i \in \mathbb{N}$ and we set $f_0\equiv 0$, with
\begin{equation}\label{3.2.2}
q^\infty_{i+1}=\left(b^\infty_{1,i+1}v_1+b^\infty_{2,i+1}v_2+c^\infty_{i+1}(\tfrac{|v|^2-3}{2})\right) \sqrt{m }.
\end{equation}
Observe that $\Gamma(f_i,f_i)\in \mathcal{N}^\perp$. By Lemma 5 in \cite{guo2010decay}, we have 
\begin{equation*}
  \left\|w\nu^{-1}\Gamma(f_i,f_i)\right\|_{L_v^\infty} \leq C\left\|wf_i\right\|^2_{L^\infty_v}.
\end{equation*}
Then by Theorem \ref{Thm-linear}, there exists a unique pair $(f_{i+1},q^\infty_{i+1})$ solving \eqref{3.2.1}, satisfying the following estimates:
\begin{equation}\label{3.2.8}
\begin{aligned}
    &|b^\infty_{1,i+1},b^\infty_{2,i+1},c^\infty_{i+1}| \leq C \left(|wR|_{L^\infty(\gamma_-)}+\left\|e^{\sigma x}\nu^{-1}wS\right\|_{L^\infty}+\left\|e^{\frac{1}{2}\sigma_0 x}wf_{i}\right\|^2_{L^\infty}\right),\\
    &\left\|e^{\sigma x}wf_{i+1}\right\|_{L^\infty}+|wf_{i+1}(0)|_{L^\infty(\gamma)}\\
    &\qquad\qquad\qquad\quad \,\ \ \leq {C}|wR|_{L^\infty(\gamma_-)} +\tfrac{{C}}{\sigma_0-\sigma}\left(\left\|e^{\sigma_0 x}\nu^{-1}wS\right\|_{L^\infty}+\left\|e^{\frac{\sigma_0 x}{2}}wf_i\right\|^2_{L^\infty}\right).
\end{aligned}
\end{equation}
We denote
\begin{equation*} 
  \delta:=|wR|_{L^\infty(\gamma_-)}+\left\|e^{\sigma_0 x}\nu^{-1}wS\right\|_{L^\infty}.
\end{equation*}
Now we assert that for all $i\geq 1$, it holds
\begin{equation}\label{3.2.3}
  \begin{aligned}
  \left|\left(b_{i, 1}^{\infty}, b_{i, 2}^{\infty}, c_{i}^{\infty}\right)\right| & \leq 2  {C} \delta, \\
  \left\|e^{\sigma x} w f_{i}\right\|_{L ^{\infty}}+\left|w f_{i}(0)\right|_{L^{\infty}(\gamma)} & \leq \tfrac{2 C \delta}{\sigma_{0}-\sigma}.
  \end{aligned}
\end{equation}
For  $i=0 $, it follows from  $f_{0} \equiv 0$ and \eqref{3.2.8} that
\begin{equation*} 
  \begin{aligned}
  \left|\left(b_{1,1}^{\infty}, b_{1,2}^{\infty}, c_{1}^{\infty}\right)\right| & \leq  {C} \delta, \\
  \left\|e^{\sigma x} w f_{1}\right\|_{L ^{\infty}}+\left|w f_{1}(0)\right|_{L^{\infty}(\gamma)} & \leq \tfrac{C \delta}{\sigma_{0}-\sigma}.
  \end{aligned}
\end{equation*}
Assume that the claim \eqref{3.2.3} holds for  $1,2, \cdots, i$, now we consider the case for  $i+1 $. It follows from \eqref{3.2.8} and \eqref{3.2.3} that
\begin{equation*} 
  \begin{aligned}
  \left|\left(b_{i+1,1}^{\infty}, b_{i+1,2}^{\infty}, c_{i+1}^{\infty}\right)\right| & \leq C\delta+C\left\|e^{\frac{1}{2} \sigma_{0} x} w f_{i}\right\|_{L ^{\infty}}^{2} \\
  & \leq C \delta\left(1+\delta\left(\tfrac{4 C}{\sigma_{0}}\right)^{2}\right) \leq \tfrac{3}{2} C \delta,
  \end{aligned}
\end{equation*}
and
\begin{equation*} 
  \begin{aligned}
  \left\|e^{\sigma x} w f_{i+1}\right\|_{L ^{\infty}}+\left|w f_{i+1}(0)\right|_{L^{\infty}(\gamma)} & \leq \tfrac{C \delta}{\sigma_{0}-\sigma}+\tfrac{C}{\sigma_{0}-\sigma}\left\|e^{\frac{1}{2} \sigma_{0} x} w f_{i}\right\|_{L ^{\infty}}^{2} \\
  & \leq \tfrac{C \delta}{\sigma_{0}-\sigma}\left(1+\delta\left(\tfrac{4 C}{\sigma_{0}}\right)^{2}\right) \leq \tfrac{\frac{3}{2} C \delta}{\sigma_{0}-\sigma},
  \end{aligned}
\end{equation*}
by choosing  $\delta$ small enough. Therefore we can conclude \eqref{3.2.3} by induction.

Finally we consider the convergence of sequence $ f_{i} $. For the difference  $f_{i+1}-f_{i} $, we apply \eqref{1.1.19} and \eqref{1.1.20} to have
\begin{equation}\label{3.2.7}
 \begin{aligned}
   &\left\|e^{\frac{1}{2} \sigma_{0} x} w\left(f_{i+1}-f_{i}\right)\right\|_{L ^{\infty}}+\left|w\left(f_{i+1}-f_{i}\right)(0)\right|_{L^{\infty}(\gamma)}  +\left|\left(b_{i+1,1}^{\infty}-b_{i, 1}^{\infty}, b_{i+1,2}^{\infty}-b_{i, 2}^{\infty}, c_{i+1}^{\infty}-c_{i}^{\infty}\right)\right| \\
   \leq &\tfrac{4 C}{\sigma_{0}}\left(\left\|e^{\sigma_{0} x} \nu^{-1} w \Gamma\left(f_{i}-f_{i-1}, f_{i}\right)\right\|_{L ^{\infty}}+\left\|e^{\sigma_{0} x} \nu^{-1} w \Gamma\left(f_{i-1}, f_{i}-f_{i-1}\right)\right\|_{L ^{\infty}}\right) \\
   \leq &\tfrac{4 C}{\sigma_{0}}\left(\left\|e^{\frac{1}{2} \sigma_{0} x} w f_{i}\right\|_{L ^{\infty}}+\left\|e^{\frac{1}{2} \sigma_{0} x} w f_{i-1}\right\|_{L ^{\infty}}\right)\left\|e^{\frac{1}{2} \sigma_{0} x} w\left(f_{i}-f_{i-1}\right)\right\|_{L ^{\infty}} \\
   \leq&  2\delta \left(\tfrac{4 C}{\sigma_{0}}\right)^{2}\left\|e^{\frac{1}{2}\sigma_0 x} w\left(f_{i}-f_{i-1}\right)\right\|_{L ^{\infty}}\\
    \leq& \tfrac{1}{2}\left\|e^{\frac{1}{2}\sigma_0 x} w\left(f_{i}-f_{i-1}\right)\right\|_{L ^{\infty}},
 \end{aligned}
\end{equation}
by taking $\delta$ sufficiently small. Hence, $\{{f}_i\}^\infty_{i=0}$ and $\{q^\infty_i\}^\infty_{i=0}$ are both Cauchy sequences, and therefore there exists $(f,q)$ solves \eqref{nonlinear-KL-eq}. The estimates \eqref{1.1.13} and \eqref{1.1.15} follows from the strong convergence and \eqref{3.2.3}. The uniqueness can also be established using \eqref{3.2.7}. Therefore we have proved Theorem \ref{Thm-nonlinear}.

\section{Application to fluid limits}\label{Sec-application}
In this section, we apply Theorem \ref{Thm-linear} to the fluid limits of the Boltzmann equation, illustrating how to derive the boundary conditions of the fluid equations by solving the Knudsen layer problem. As shown below, the derivation of the boundary conditions essencially relies on the symmetric properties of the Knudsen layer problem. Specifically, we consider the following linearized Boltzmann in half-space with incompressible scaling:
\begin{equation}\label{linear-BE}
  \varepsilon \partial_t f^\varepsilon+v\cdot \nabla_xf^\varepsilon+\frac{1}{\varepsilon}\mathcal{L}f^\varepsilon=0, \quad (x,v)\in \mathbb{R}^3_+\times \mathbb{R}^3.
\end{equation}
Here $\varepsilon>0$ represents the Knudsen number, which denotes the ratio of the mean free path to the macroscopic length scale. We impose \eqref{linear-BE} with the Maxwell reflection boundary condition
\begin{equation}\label{BE-BC}
  (L^R-\alpha L^D)f^\varepsilon=0,
\end{equation}
where $L^R$ and $L^D$ are defined in Section \ref{Se1.1}. We assume $\alpha=O(1), \alpha\in(0,1]$. Therefore, the case of specular or almost specular reflection boundary conditions is excluded and the case of diffuse reflection is included. We seek the solution to \eqref{linear-BE} in the following form
\begin{equation}\label{expansion}
\begin{aligned}
    f^\varepsilon=f_0+&\varepsilon f_1+\cdots\\
    +&\varepsilon f^{bb}_1+\cdots,
\end{aligned}
\end{equation}
where $f_1^{bb}=f_1^{bb}(t,x_1,x_2,\xi,v)$ with $\xi=\frac{x_3}{\varepsilon}$ is the Knudsen layer term. Substituting \eqref{expansion}  into \eqref{linear-BE} and rearranging by the order of $\eps$, we obtain (see \cite{takata2012asymptotic} for details of the computation)
\begin{equation}\label{4.1.3}
  f_0=\left(\rho_0+u_0 \cdot v+\theta_0(\tfrac{|v|^2-3}{2})\right)\sqrt{m},
\end{equation}
with $\nabla_x\cdot u_0=0$ and $\nabla_x(\rho_0+\theta_0)=0$. Moreover, $(u_0,\theta_0)$ satisfy the linearized incompressible Navier-Stokes-Fourier system. The first order can be expressed as
\begin{equation}\label{4.1.4}
  f_1=\left(\rho_1+u_1 \cdot v+\theta_1(\tfrac{|v|^2-3}{2})\right)\sqrt{m}-\tfrac{1}{2}\hat{A}:  \sigma(u_0)-\hat{B}:\nabla_x \theta_0,
\end{equation}
 where $\hat{A}$ and $\hat{B}$ are defined in \eqref{hatAB} and $\sigma(u_0)_{ij}=\partial_iu_{0,j}+\partial_ju_{0,i}-\tfrac{2}{3}\nabla_x\cdot u_0$. Moreover, the boundary layer term $f_1^{bb}$ satisfies
\begin{equation}\label{KL-fbb}
\begin{cases}
    v_3\partial_\xi f^{bb}_1+\mathcal{L}f^{bb}_1=0, \quad (\xi>0),\\
    (L^R-\alpha L^D)f^{bb}_1=-(L^R-\alpha L^D)f_1, \quad (\xi=0,v_3>0),\\
    f_1^{bb}\to 0 \quad \text{as } \xi\to \infty.
\end{cases}
\end{equation}
From now on, we focus on deriving the boundary conditions of the fluid quantities above.

$O(1)$: Substituting \eqref{4.1.3} into the Maxwell reflection boundary condition \eqref{BE-BC}, we obtain
\begin{equation*}
  (L^R-\alpha L^D)f_0=0.
\end{equation*}
This equation further leads to
\begin{equation}\label{4.1.7}
  u_0=0,\quad \theta_0=0,\quad \text{on} \quad \! \partial \mathbb{R}^3_+.
\end{equation}

$O(\eps)$: Since the boundary is impermeable, we deduce that $\int_{\R^3}v_3 (f_1+f^{bb}_1) d v=0$ at $x_3=0$. Moreover, by the conservation law in  \eqref{KL-fbb}, it follows that $\int_{\R^3}v_3 f^{bb}_1 d v \equiv 0$. Therefore 
\begin{equation}\label{u13}
  u_{1,3}=0,\quad \text{at}\quad\! x_3=0
\end{equation}
Next, substituting \eqref{4.1.4} into \eqref{BE-BC}, one has
\begin{equation}\label{4.1.9}
  \begin{cases}
      v_3\partial_\xi f^{bb}_1+\mathcal{L}f^{bb}_1=0, \quad (\xi>0),\\
      \begin{aligned}
        (L^R-\alpha L^D)f^{bb}_1=&-\sum_{i=1}^{2}\alpha u_{1,i}v_i\sqrt{m(v)}-\alpha\theta_1(\tfrac{|v|^2}{2}-2)\sqrt{m(v)}\\
        &+(2-\alpha)\sum_{i=1}^{2}\hat{A}_{i3}\left(\partial_iu_{0,3}+\partial_3u_{0,i}\right)+(2-\alpha)\hat{B}_3\partial_3\theta_0, \quad (\xi=0,v_3>0),
      \end{aligned}\\
      f^{bb}_1 \to 0, \quad \text{as }\xi \to \infty,
  \end{cases}
  \end{equation}
where we use the the fact that $\nabla_x \cdot u_0=0$ and \eqref{4.1.7}. It is important to note that $u_0$ and $\theta_0$ are now determined. By Theorem \ref{Thm-linear}, there exists a unique pair $(u_{1,1}, u_{1,2},\theta_1)$ such that \eqref{4.1.9} has a unique solution. To obtain explicit expressions of $u_{1,i}(i=1,2)$ and $\theta_1$, we use the symmetric properties of the linearized Boltzmann operator $\L$ and the boundary operator $L^R-\alpha L^D$. Specifically, let $f(v)=f(|v|,v_3)$ be a function that is radial in $v_1$ and $v_2$, we now consider the symmetric tensor fields
  $$G_i( v )=\mathcal{L} v _if(| v |, v_3) ,\quad\!i=1,2, \quad  \text{and}\quad \! G( v )=\mathcal{L} f(| v |, v_3) ,$$
Based on the isotropic property of $\mathcal{L}$, both $G_i$ and $G$ are tensor fields axially symmetric with respect to axis $(0,0,-1)$(see \cite{sone2002kinetic} Appendix B and \cite{sone2007molecular} Appendix A.2). There exists a unique $g(|v|,v_3)$ and $g_0(|v|,v_3)$, both radial in $v_1$ and $v_2$, such that :
\begin{equation*}
  G_i( v )=v _ig_0(| v |, v_3),\quad\!i=1,2, \quad \text{and} \quad \!  G( v )=g(| v |, v_3).
  \end{equation*}
  Therefore, we define the operator $\mathcal{L}^S$ as follows:
  \begin{equation}\label{LS}
    \mathcal{L}[ v _if(| v |, v_3)]=  v _ig_0(| v |, v_3):=v _i\mathcal{L}^Sf(| v |, v_3),i=1,2.
  \end{equation}
  Based on the property above, the operator $v_3  \partial_\xi +\mathcal{L}$ conserves the symmetry with respect to $v _i$ for $i=1,2$. Specifically,
   \begin{equation*} 
    \left(v_3\partial_\xi+\mathcal{L}\right)f(| v |, v_3)=g(| v |, v_3),
  \end{equation*}
   and
   \begin{equation*} 
     \left(v_3 {\partial_\xi}+\mathcal{L}\right)v _if_i(| v |, v_3)=v _ig_i(| v |, v_3),\quad i=1,2.
   \end{equation*}
  Furthermore, a direct calculation shows that $L^R-\alpha L^D$ satisfies the same property. We apply this property to the Knudsen-Layer problem. According to the inhomogeneous terms in \eqref{4.1.9}, we seek the solution of \eqref{4.1.9} in the following form:
  \begin{equation}\label{4.1.12}
  \begin{aligned}
  f^{bb}_1(t,x_1,x_2, \xi ,v)=&(\partial_iu_{0,3}+\partial_3u_{0,i})_Bv_i \phi_u(\xi,| v |, v_3)\\
  &+(\partial_3 \theta_0)_B\phi_\theta(\xi,| v |, v_3),
  \end{aligned}
  \end{equation}
  and
  \begin{equation}\label{4.1.13}
  \begin{aligned}
    &(u_{1,i})_B=c_u \left(\partial_3u_{0,i}+\partial_iu_{0,3}\right)_B ,\
    &(\theta_1)_B=c_\theta (\partial_3\theta_0)_B.
  \end{aligned}
\end{equation}
Here the subscript $B$ denotes the value on the boundary, while $c_u$ and $c_\theta$ are undetermined constants to be determined together with the solutions $\phi_u$ and $\phi_\theta$, as stated in Theorem \ref{Thm-linear}. The assumption that $\phi_u$ and $\phi_\theta$ depend only on $\xi,v_3$ and $| v |$ is shown to be consistent.

Substituting \eqref{4.1.12} and \eqref{4.1.13} into the Knudsen-Layer problem \eqref{4.1.9} and using \eqref{LS}, we obtain the problem for $\phi_u$ and $\phi_\theta$. That is,
  \begin{itemize}
    \item Problem for $\phi_u$
    \begin{equation}\label{4.1.14}
      \left\{
    \begin{aligned}
    &(v_3 \partial_\xi +\mathcal{L}^S)\phi_u(\xi,v)=0, \qquad(\xi>0),\\
    &\phi_u(0,v)|_{v_3>0}=(1-\alpha)L\gamma_+\phi_u(0,v)-\alpha c_u\sqrt{m( v )}+(2-\alpha)v_3\alpha(| v |)\sqrt{m( v )}, \\
    &\lim_{\xi\to\infty}\phi_u = 0.
    \end{aligned}\right.
    \end{equation}
   \item Problem for $\phi_\theta$
    \begin{equation}\label{4.1.15}
      \left\{
    \begin{aligned}
    &(v_3 \partial_\xi +\mathcal{L})\phi_\theta(\xi,v)=0, \qquad(\xi>0),\\
    &(L^R-\alpha L^D)\phi_\theta(0,v)=-\alpha(\tfrac{| v |^2}{2}-2) c_\theta\sqrt{m( v )}+(2-\alpha)v_3 \tfrac{| v |^2-5}{2}  b(| v |)\sqrt{m( v )},\\
    &\lim_{\xi\to \infty}\phi_\theta =0.
    \end{aligned}\right.\end{equation}
  \end{itemize}
In the following, we will demonstrate that both \eqref{4.1.14} and \eqref{4.1.15} are solvable. To solve \eqref{4.1.15}, we present the following lemma:
\begin{lemma}\label{Lm4.1}
 Assume that $\alpha=O(1)$ and $0<\alpha\leq  1$. There exists a unique pair $(c_\theta,\phi_\theta)$ solves \eqref{4.1.15}, and satisfies the following estimates.
  \begin{equation}\label{4.0.2}
    |c_\theta|\leq C , \quad \text{and} \quad \!\|e^{\sigma x}\nu^{-1} w \phi_\theta\|_{L^\infty_{x,v}} \leq C ,
  \end{equation}
  for any sufficiently small $\sigma>0$.
\end{lemma}
\begin{proof}
  First, we consider the following system
  \begin{equation}
    \left\{
    \begin{aligned}\label{4.0.1}
    &(v _3 \partial_\xi+\mathcal{L}) \phi(\xi ,v )=0, \qquad(\xi>0),\\
    &(L^R-\alpha L^D)\phi( \xi,v )=-\alpha(\tfrac{| v |^2}{2}-2)C\sqrt{m( v )}+r(| v |, v _3)\qquad(\xi=0, v _3>0),
    \end{aligned}\right.
    \end{equation}
  where $r(|v|,v_3)$ is a function radially symmetric in $v_1,v_2$, satisfying the conditions \eqref{1.1.8} and \eqref{1.1.25}. We will show that the solution of \eqref{4.0.1} is also radially symmetric in $v_1,v_2$, i.e., $\phi(\xi,v)=\phi(\xi,|v|, v_3)$. To prove this, we introduce the rotation of coordinates $ v _1$ and $ v _2$ as follows:
  $$R v =
  \begin{pmatrix}
   \sin \theta  & \cos \theta  &0 \\
   -\cos \theta  & \sin \theta  & 0\\
   0 & 0 & 1
  \end{pmatrix}
  \begin{pmatrix}
   v  _1\\
   v _2 \\
   v _3
  \end{pmatrix}
  =\begin{pmatrix}
   \zeta _1\\
  \zeta_2 \\
  \zeta_3
  \end{pmatrix},
  $$
  for $\theta \in[0,2\pi]$. A simple calculation yields that, $\phi(\xi,R v)$ satisfies the same system as $\phi(\xi,v)$, i.e.,
  \begin{equation*}
  \left\{
  \begin{aligned}
  & (v _3 \partial_\xi+L) \phi(\xi,R v )=0, \qquad(\xi>0),\\
  &(L^R-\alpha L^D)\phi(\xi,R v )=-\alpha(\tfrac{| v |^2}{2}-2)C\sqrt{m( v )}+r(| v |, v _3)\qquad(\xi=0, v _3>0),
  \end{aligned}\right.
  \end{equation*}
  where we use the fact that $\mathcal{L}(\phi(R v ))=[\mathcal{L}\phi](R v )$, see \cite{sone2007molecular} for example. Since there is no uniqueness theorem for the system \eqref{4.0.1}, we cannot conclude that $\phi(\xi,R v)=\phi(\xi,v)$. Thanks to Lemma 2.1.1 in \cite{coron1988classification}, we consider function $\Phi(\xi, v )=\phi(\xi, v )-\phi(\xi,R v )$. It is straightforward to verify that $\Phi$ satisfies the condition $\int_{ v _3<0} v _3\Phi(0, v )\sqrt{m}d v =0$, and the system for $\Phi$ is
  \begin{equation*} 
  \left\{
  \begin{aligned}
  &  (v _3 \partial_ \xi +\mathcal{L}) \Phi( \xi , v )=0,\\
  &(L^R-\alpha L^D)\Phi(0, v )=0,\\
  &\int_{ v _3<0} v _3\Phi(0, v )\sqrt{m}dv =0.
  \end{aligned}\right.
  \end{equation*}
  It immediately follows that $\Phi(\xi, v ) \equiv 0$ by the uniqueness result in Lemma 2.1.1 of \cite{golse1988boundary}. This indicates that $\phi$ is rotation-invariant with respect to $ v _1$ and $ v _2$, and $\phi$ is a $| v |, v _3-$dependence function, i.e.,
  \begin{equation}\label{4.0.4}
    \phi(\xi, v )=\phi(\xi,| v |, v _3).
  \end{equation}
  Therefore, the oddness implies that
  \begin{equation}\label{4.0.5}
    \mathbf{P}\phi=\left(a+c (\tfrac{|v|^2-3}{2} )\right)\sqrt{m(v)}.
  \end{equation}
  For $0<\alpha<1$, Theorem \ref{Thm-linear}, \eqref{4.0.4} and \eqref{4.0.5} imply that there exists a unique pair $(\phi_\theta,c_\theta)$ solves \eqref{4.1.15}. Moreover, the estimates in \eqref{4.0.2} follow from Theorem \ref{Thm-linear}. On the other hand, for $\alpha=1$, by using Theorem 1.5 in \cite{huang2022boundary}, we can also obtain the existence and \eqref{4.0.2}. The proof of Lemma \ref{Lm4.1} is therefore complete.
\end{proof}

To solve problem \eqref{4.1.14}, we first list some properties of $\mathcal{L}^S$.
\begin{proposition}\label{Pr4.1}(1) $\mathcal{L}^S$ is a self-adjoint operator on $L^2(d\mu)$, where $d\mu=v_r^3dv_rdv_3$ and $v_r:=\sqrt{v_1^2+v_2^2}\in[0,\infty)$:
   $$\int\mathcal{L}^S f(v_r ,v_3)g(v_r ,v_3) d\mu =\int L^Sg(v_r ,v_3)f(v_r ,v_3) d\mu ,$$
for $f(v_r,v_3),g(v_r,v_3)\in L^2(d\mu)$.

(2)The null space of $\L^S$ is $Null(\mathcal{L}^S)=span\{\sqrt{m(v)}\}$.

(3)Coercivity of $\L^S$: Let $\varphi(v_r,v_3)\in Null(\mathcal{L}^S)^\perp_{ d\mu }$, i.e. $\int\varphi(v_r,v_3)\sqrt{m(v)}d\mu=0$, then
$$(\mathcal{L}^S\varphi(v_r,v_3),\varphi(v_r,v_3))_{ d\mu } \geq c_0\int \nu \varphi(v_r,v_3)^2 d\mu.$$

(4)$\mathcal{L}^S=\nu(|v|)-K^S$, where $\nu(|v|)$ is the collision frequency, defined in \eqref{coll-freq} and $K^S$ is a compact operator on $L^2(d\mu)$.
\end{proposition}
\begin{proof}
  (1) By the definition of $\L^S$ in \eqref{LS}, we have for any $c_1,c_2\in \R$
  \begin{equation}
    \mathcal{L}\left((c_1v_1+c_2v_2)\varphi(v_r,v_3)\right)=(c_1v_1+c_2v_2)\mathcal{L}^S\left(\varphi(v_r,v_3)\right).
  \end{equation}
   Let $v_1=v_r \sin \theta,v_2=v_r\cos\theta$ for $\theta\in[0,2\pi]$. We then have for any $c_1,c_2\in \R$,
  \begin{equation*} 
   \begin{aligned}
     \int_{0}^{2\pi}\left(c_1\sin \theta+c_2\cos \theta\right)^2d\theta \int_{\mathbb{R}^2_+}\mathcal{L}^S\varphi_1(|v|,v_3) \varphi_2(|v|,v_3) d\mu \\
     =\int_{\mathbb{R}^3}\mathcal{L}\left((c_1v_1+c_2v_2)\varphi_1(v_r,v_3)\right)\left((c_1v_1+c_2v_2)\varphi_2(v_r,v_3)\right)dv\\
     =\int_{\mathbb{R}^3}\mathcal{L}\left((c_1v_1+c_2v_2)\varphi_2(v_r,v_3)\right)\left((c_1v_1+c_2v_2)\varphi_1(v_r,v_3)\right)dv\\
     =   \int_{0}^{2\pi}\left(c_1\sin \theta+c_2\cos \theta\right)^2d\theta \int_{\mathbb{R}^2_+}\mathcal{L}^S\varphi_2(|v|,v_3) \varphi_1(|v|,v_3) d\mu,
   \end{aligned}
  \end{equation*}
  where we use the fact that $\mathcal{L}$ is a self-adjoint operator on $L^2(\mathbb{R}^3)$. Here, we choose $c_1^2+c_2^2 \ne 0$, and therefore
    $\int_{0}^{2\pi}\left(c_1\sin \theta+c_2\cos \theta\right)^2d\theta\ne 0.$
Hence,
  \begin{equation*} 
       \int_{\mathbb{R}^2_+}\mathcal{L}^S\varphi_1(|v|,v_3) \varphi_2(|v|,v_3) d\mu 
     =  \int_{\mathbb{R}^2_+}\mathcal{L}^S\varphi_2(|v|,v_3) \varphi_1(|v|,v_3) d\mu.
   \end{equation*}
This implies that $\mathcal{L}^S$ is a self-adjoint operator on $L^2(d\mu)$.
  
  (2)Noting that
  \begin{equation*}
    \mathcal{L}^S\varphi(v_3,|v|)=0 \quad \text{if and only if} \quad\!\mathcal{L}\left((c_1v_1+c_2v_2)\varphi(v_3,|v|)\right)=0,
  \end{equation*}
  for any $c_1,c_2\in \R$. Then we have
  \begin{equation*}
    (c_1v_1+c_2v_2)\varphi(v_r,v_3) \in \text{span}\{\sqrt{m},v\sqrt{m},\tfrac{|v|^2-3}{2}\sqrt{m}\}.
  \end{equation*}
  The oddness implies that
  \begin{equation*}
    \int_{\mathbb{R}^3}(c_1v_1+c_2v_2)\varphi(v_r,v_3)\sqrt{m}dv= 0,
  \end{equation*}
  \begin{equation*}
    \int_{\mathbb{R}^3}(c_1v_1+c_2v_2)\varphi(v_r,v_3)\tfrac{|v|^2-3}{2}\sqrt{m}dv= 0,
  \end{equation*}
  and
  \begin{equation*}
    \int_{\mathbb{R}^3}(c_1v_1+c_2v_2)\varphi(v_r,v_3)v_3\sqrt{m}dv= 0.
  \end{equation*}
  We choose $c_2=0$ to obtain
  \begin{equation*}
    c_1v_1\varphi(v_r,v_3)=Cv_1\sqrt{m},
  \end{equation*}
  which yields that
  \begin{equation*}
    \varphi(v_r,v_3)\in \text{span}\{\sqrt{m}\}.
  \end{equation*}
  
  (3)Suppose $\int\varphi(v_r,v_3)\sqrt{m}d\mu=0$. Then we have
  \begin{equation*}
    \int_{\mathbb{R}^3}(c_1v_1)\varphi(v_r,v_3)v_1\sqrt{m(v)}dv=0 \quad \text{for}\quad \! c_1\neq 0.
  \end{equation*}
  Combining this with the oddness of $c_1v_1\varphi(v_r,v_3)$, we conclude that $c_1v_1\varphi(v_r,v_3)\in \mathcal{N}^\perp$. The coercivity of $\L$ implies that
  \begin{equation*}
    \int_{\mathbb{R}^3}\mathcal{L}\left(c_1v_1\varphi(v_r,v_3)\right)c_1v_1\varphi(v_r,v_3)dv \geq c_0\|c_1v_1\varphi(v_r,v_3)\|^2_{\nu}.
  \end{equation*}
  Consequently,
  \begin{equation*}
    (\mathcal{L}^S \varphi,\varphi)_{d\mu}\geq c_0\|\varphi\sqrt{\nu}\| ^2_{L^2(d\mu)}.
  \end{equation*}
  Thus, the coercivity of $\L^S$ is established.

  (4)We observe that for any $c_1,c_2\in \mathbb{R}$,
  \begin{equation*}
   \begin{aligned}
     \mathcal{L}\left((c_1v_1+c_2v_2)\varphi\right)=&\nu(|v|)(c_1v_1+c_2v_2)\varphi-K\left((c_1v_1+c_2v_2)\varphi\right)\\
     =&(c_1v_1+c_2v_2)\mathcal{L}^S\varphi(v_r,v_3).
   \end{aligned}
  \end{equation*}
  We define the operator $K^S: L^2(d \mu)\to L^2(d \mu)$ as
  \begin{equation*}
   K^S\varphi(v_r,v_3)=\nu(|v|)\varphi(v_r,v_3)-\mathcal{L}^S\varphi(v_r,v_3),
  \end{equation*}
  Then $K^S$ satisfies that
  \begin{equation*}
    (c_1v_1+c_2v_2)K^S\varphi(v_r,v_3)=K\big((c_1v_1+c_2v_2)\varphi(v_r,v_3)\big).
  \end{equation*}
   The compactness of $K^S$ on $L^2(d\mu)$ follows from the compactness of $K$ on $L^2(\mathbb{R}^3)$.
\end{proof}

\begin{lemma}\label{Lm4.3}
  Assume that $\alpha= O(1)$ and $0 <\alpha\leq 1$. There exists a unique pair $(\phi_u,c_u)$ solves \eqref{4.1.14}. Moreover, the following estimates hold:
  \begin{equation*}
    \|e^{\sigma \xi} \nu^{-1}w v_i\phi_u\|_{L^\infty} \leq C, \quad (i=1,2) \quad \text{and} \quad\!|c_u| \leq C,
  \end{equation*}
for any sufficiently small $\sigma$.
\end{lemma}
  \begin{proof}
    Proposition \ref{Pr4.1} and \eqref{LS} imply that the operator $\L^S$ has similar properties to the linearized Boltzmann operator $\L$. Then by similar arguments as in \cite{bardos1986milne} and \cite{coron1988classification}, we can prove there exists a unique $\phi_u,c_u$ solves \eqref{4.1.14} such that $\phi_u\in L^\infty\left( d\xi ;L^2(|v_3|d\mu)\right)$. It is straightforward to see that $(\phi_u, c_u)$ satisfies
  \begin{equation}
    \begin{cases}
      v_3\partial_\xi(v_i\phi_u)+\mathcal{L}(v_i\phi_u)=0,\\
      (L^R-\alpha L^D)v_i\phi_u(0,v)=-\alpha c_uv_i\sqrt{m}+(2-\alpha)v_3v_i\alpha(|v|)\sqrt{m},\\
      \lim_{\xi\to\infty}v_i\phi_u = 0,
    \end{cases}
  \end{equation}
  for $i=1,2$. For $0<\alpha<1$, we apply Theorem \ref{Thm-linear} to obtain that
  \begin{equation}\label{4.2.3}
    \begin{array}{l}
      \|e^{\sigma \xi} \nu^{-1}w v_i\phi_u\|_{L^\infty} \leq C ,\\
      |C_u|\leq C .
    \end{array}
  \end{equation}
  For the case of $\alpha=1$, we can also obtain \eqref{4.2.3} by using Theorem 1.5 in \cite{huang2022boundary}. Therefore the prove of Lemma \ref{Lm4.3} is complete.
  \end{proof}

Based on Lemmas \ref{Lm4.1} and \ref{Lm4.3}, we have constructed a solution to \eqref{4.1.9}. Furthermore, we know that \eqref{4.1.12} and \eqref{4.1.13} is the unique solution to \eqref{4.1.9} due to Theorem \ref{Thm-linear} and Theorem 1.5 in \cite{huang2022boundary}.

From \eqref{u13} and \eqref{4.1.13}, we can immediately write down the slip boundary conditions for the first order in the following form:
\begin{equation*}
 \begin{cases}
  u_{1,i}=c_u\left(\partial_3u_{0,i}+\partial_i u_{0,3}\right), \quad i=1,2,\\
u_{1,3}=0,\\
\theta_1=c_\theta \partial_3\theta_0,
 \end{cases}
\end{equation*}
where $c_u$ and $c_\theta$ are the so-called slip coefficients. For hard-sphere molecules, the slip coefficients are uniquely determined by the accommodation coefficient $\alpha$. Additionally, the specific values of these slip coefficients have been summarized by the authors in \cite{sone2002kinetic,sone2007molecular,takata2012asymptotic,aoki2017slip}.

In summary, the boundary conditions of the fluid equations can be derived through the analysis of the Knudsen layer, which heavily relies on the symmetric properties of $\L$ and the boundary operator $L^R-\alpha L^D$. Moreover, since the nonlinear term $\Gamma$ possesses similar symmetric properties to $\L$, we can apply this method to the hydrodynamic limit of the nonlinear Boltzmann equation.

\bibliography{reff} 

\begin{thebibliography}{10}

\bibitem{aoki2017slip}
K.~Aoki, C.~Baranger, M.~Hattori, S.~Kosuge, G.~Martal\`o, J.~Mathiaud, and L.~Mieussens.
\newblock Slip boundary conditions for the compressible {N}avier-{S}tokes equations.
\newblock {\em J. Stat. Phys.}, 169(4):744--781, 2017.

\bibitem{bardos1986milne}
C.~Bardos, R.~E. Caflisch, and B.~Nicolaenko.
\newblock The {M}ilne and {K}ramers problems for the {B}oltzmann equation of a hard sphere gas.
\newblock {\em Comm. Pure Appl. Math.}, 39(3):323--352, 1986.

\bibitem{Niclas-Golse-2021-ARMA}
N.~Bernhoff and F.~Golse.
\newblock On the boundary layer equations with phase transition in the kinetic theory of gases.
\newblock {\em Arch. Ration. Mech. Anal.}, 240(1):51--98, 2021.

\bibitem{caflisch1980fluid}
R.~E. Caflisch.
\newblock The fluid dynamic limit of the nonlinear {B}oltzmann equation.
\newblock {\em Comm. Pure Appl. Math.}, 33(5):651--666, 1980.

\bibitem{Chen-Liu-Yang-2004-AA}
C.-C. Chen, T.-P. Liu, and T.~Yang.
\newblock Existence of boundary layer solutions to the {B}oltzmann equation.
\newblock {\em Anal. Appl. (Singap.)}, 2(4):337--363, 2004.

\bibitem{coron1988classification}
F.~Coron, F.~Golse, and C.~Sulem.
\newblock A classification of well-posed kinetic layer problems.
\newblock {\em Comm. Pure Appl. Math.}, 41(4):409--435, 1988.

\bibitem{deng2008pointwise}
S.~Deng, W.~Wang, and S.-H. Yu.
\newblock Pointwise convergence to {K}nudsen layers of the {B}oltzmann equation.
\newblock {\em Comm. Math. Phys.}, 281(2):287--347, 2008.

\bibitem{desvillettes1994remark}
L.~Desvillettes and F.~Golse.
\newblock A remark concerning the {C}hapman-{E}nskog asymptotics.
\newblock In {\em Advances in kinetic theory and computing}, volume~22 of {\em Ser. Adv. Math. Appl. Sci.}, pages 191--203. World Sci. Publ., River Edge, NJ, 1994.

\bibitem{esposito2013non}
R.~Esposito, Y.~Guo, C.~Kim, and R.~Marra.
\newblock Non-isothermal boundary in the {B}oltzmann theory and {F}ourier law.
\newblock {\em Comm. Math. Phys.}, 323(1):177--239, 2013.

\bibitem{golse1988boundary}
F.~Golse, B.~Perthame, and C.~Sulem.
\newblock On a boundary layer problem for the nonlinear {B}oltzmann equation.
\newblock {\em Arch. Rational Mech. Anal.}, 103(1):81--96, 1988.

\bibitem{golse1989stationary}
F.~Golse and F.~Poupaud.
\newblock Stationary solutions of the linearized {B}oltzmann equation in a half-space.
\newblock {\em Math. Methods Appl. Sci.}, 11(4):483--502, 1989.

\bibitem{guo2010decay}
Y.~Guo.
\newblock Decay and continuity of the {B}oltzmann equation in bounded domains.
\newblock {\em Arch. Ration. Mech. Anal.}, 197(3):713--809, 2010.

\bibitem{takata2012asymptotic}
M.~Hattori and S.~Takata.
\newblock Asymptotic theory for the time-dependent behavior of a slightly rarefied gas over a smooth solid boundary.
\newblock {\em J. Stat. Phys.}, 147:1182--1215, 2012.

\bibitem{huang2022boundary}
F.~Huang and Y.~Wang.
\newblock Boundary layer solution of the {B}oltzmann equation for diffusive reflection boundary conditions in half-space.
\newblock {\em SIAM J. Math. Anal.}, 54(3):3480--3534, 2022.

\bibitem{huang2023boundary}
F.-M. Huang, Z.-H. Jiang, and Y.~Wang.
\newblock Boundary layer solution of the {B}oltzmann equation for specular boundary condition.
\newblock {\em Acta Math. Appl. Sin. Engl. Ser.}, 39(1):65--94, 2023.

\bibitem{jiang2021compressiblea}
N.~Jiang, Y.-L. Luo, and S.~Tang.
\newblock Compressible {E}uler limit from {B}oltzmann equation with {M}axwell reflection boundary condition in half-space.
\newblock {\em arXiv:2101.11199}.

\bibitem{jiang2021compressibleb}
N.~Jiang, Y.-L. Luo, and S.~Tang.
\newblock Compressible {E}uler limit from {B}oltzmann equation with complete diffusive boundary condition in half-space.
\newblock {\em Trans. Amer. Math. Soc.}, 377(8):5323--5359, 2024.

\bibitem{jiang2024knudsenboundarylayerequations}
N.~Jiang, Y.-L. Luo, and Y.~Wu.
\newblock Knudsen boundary layer equations for full ranges of cutoff collision kernels: Maxwell reflection boundary with all accommodation coefficients in [0,1].
\newblock {\em arXiv:2407.02852}.

\bibitem{Masmoudi-Jiang-2017-CPAM}
N.~Jiang and N.~Masmoudi.
\newblock Boundary layers and incompressible {N}avier-{S}tokes-{F}ourier limit of the {B}oltzmann equation in bounded domain {I}.
\newblock {\em Comm. Pure Appl. Math.}, 70(1):90--171, 2017.

\bibitem{Liu-Yu-2013-ARMA}
T.-P. Liu and S.-H. Yu.
\newblock Invariant manifolds for steady {B}oltzmann flows and applications.
\newblock {\em Arch. Ration. Mech. Anal.}, 209(3):869--997, 2013.

\bibitem{Shota-2022-JDE}
S.~Sakamoto, M.~Suzuki, and K.~Z. Zhang.
\newblock Boundary layers of the {B}oltzmann equation in three-dimensional half-space.
\newblock {\em J. Differential Equations}, 314:446--472, 2022.

\bibitem{Sammartino-Caflisch-1988-CMP1}
M.~Sammartino and R.~E. Caflisch.
\newblock Zero viscosity limit for analytic solutions, of the {N}avier-{S}tokes equation on a half-space. {I}. {E}xistence for {E}uler and {P}randtl equations.
\newblock {\em Comm. Math. Phys.}, 192(2):433--461, 1998.

\bibitem{Sammartino-Caflisch-1988-CMP2}
M.~Sammartino and R.~E. Caflisch.
\newblock Zero viscosity limit for analytic solutions of the {N}avier-{S}tokes equation on a half-space. {II}. {C}onstruction of the {N}avier-{S}tokes solution.
\newblock {\em Comm. Math. Phys.}, 192(2):463--491, 1998.

\bibitem{sone2002kinetic}
Y.~Sone.
\newblock {\em Kinetic theory and fluid dynamics}.
\newblock Modeling and Simulation in Science, Engineering and Technology. Birkh\"{a}user Boston, Inc., Boston, MA, 2002.

\bibitem{sone2007molecular}
Y.~Sone.
\newblock {\em Molecular gas dynamics}.
\newblock Modeling and Simulation in Science, Engineering and Technology. Birkh\"{a}user Boston, Inc., Boston, MA, 2007.
\newblock Theory, techniques, and applications.

\bibitem{ukai2003nonlinear}
S.~Ukai, T.~Yang, and S.-H. Yu.
\newblock Nonlinear boundary layers of the {B}oltzmann equation. {I}. {E}xistence.
\newblock {\em Comm. Math. Phys.}, 236(3):373--393, 2003.

\bibitem{ukai2004nonlinear}
S.~Ukai, T.~Yang, and S.-H. Yu.
\newblock Nonlinear stability of boundary layers of the {B}oltzmann equation. {I}. {T}he case {$\mathcal{M}^\infty<-1$}.
\newblock {\em Comm. Math. Phys.}, 244(1):99--109, 2004.

\bibitem{Wang-Yang-Yang-2006-JMP}
W.~Wang, T.~Yang, and X.~Yang.
\newblock Nonlinear stability of boundary layers of the {B}oltzmann equation for cutoff hard potentials.
\newblock {\em J. Math. Phys.}, 47(8):083301, 15, 2006.

\bibitem{Wang-Yang-Yang-2007-JMP}
W.~Wang, T.~Yang, and X.~Yang.
\newblock Existence of boundary layers to the {B}oltzmann equation with cutoff soft potentials.
\newblock {\em J. Math. Phys.}, 48(7):073304, 21, 2007.

\bibitem{Yang-2011-JSP}
X.~Yang.
\newblock The solutions for the boundary layer problem of {B}oltzmann equation in a half-space.
\newblock {\em J. Stat. Phys.}, 143(1):168--196, 2011.

\end{thebibliography}
\bibliographystyle{abbrv} 

\end{document}